\documentclass[a4paper]{article}

\usepackage[english]{babel}
\usepackage[utf8x]{inputenc}
\usepackage[T1]{fontenc}
\usepackage{lipsum}
\usepackage{blindtext}
\usepackage{graphicx,amssymb,amsmath,textcomp}
\usepackage{newpxtext} % ``New PX" font with original code \usepackage{newpxtext,newpxmath}

\usepackage[a4paper,top=3cm,bottom=2cm,left=3cm,right=3cm,marginparwidth=1.75cm]{geometry}
\usepackage{tikz}
\usetikzlibrary{matrix,arrows,positioning,calc}
\usepackage{verbatim}

\usepackage{amsmath,amsthm}
\usepackage{faktor}
\usepackage{xfrac} 
\usepackage{graphicx}
\usepackage[colorinlistoftodos]{todonotes}
\usepackage[colorlinks=true, allcolors=blue]{hyperref}
\usepackage{mathtools}
\usepackage{tikz-cd}
\usepackage{cleveref}

\newcommand{\Hom}{{\rm Hom}}
\newcommand{\vdim}{{\rm vdim \ }}
\newcommand{\Ima}{{\rm Im}}
\newcommand{\coker}{{\rm coker  }}
\newcommand{\rank}{{\rm rank \ }}
\newcommand{\cspec}{{\rm spec}}
\newcommand{\spec}{{\rm Spec}}
\newcommand{\Span}{{\rm Span}}

\newcommand{\N}{{\mathbb N}}
\newcommand{\Z}{{\mathbb Z}}
\newcommand{\R}{{\mathbb R}}

\newcommand{\K}{{\mathbb K}}

\newcommand{\dR}{{d_{dR}}}
\title{\textbf{Shifted Contact Structures  and Their Local Theory }}
\author{Kadri İlker Berktav \footnote{University of Zurich, Institute of Mathematics, Zurich, Switzerland; e-mail: kadriilker.berktav@math.uzh.ch } \\ %Zurich, Switzerland; e-mail(s): kadiilker.berktav@math.uzh.ch, ilkerberktav@gmail.com \textit{Department of Mathematics, Middle East Technical University}, \\ \textit{06800 Ankara, Turkey}
} 

\date{\vspace{-5ex}}
\begin{document}
%%%%%%%%%%%%%%%%%%%% Text italic %%%%%%%%%%%%%%%%%%%%%%%%%%%%
\theoremstyle{plain}
\newtheorem{theorem}{Theorem}[section] %[section] or [subsection]
\newtheorem{lemma}[theorem]{Lemma} 
%NOT: araya [theorem] yazdigimiz zaman Lemma/propositon etc.'lari kendi icinde siralamak yerine 2. kod ile verilen genel numaralandirmayi (section veya subsection) takip eder.
\newtheorem{proposition}[theorem]{Proposition} %[section] or [subsection]
\newtheorem{corollary}[theorem]{Corollary}%[section]
\newtheorem{conjecture}[theorem]{Conjecture}%[section]
\newtheorem{claim}[theorem]{Claim}%[section]
%%%%%%%%%%%%%%%%%%%% Text roman %%%%%%%%%%%%%%%%%%%%%%%%%%%%%
\theoremstyle{definition}
\newtheorem{notations}[theorem]{Notations}
\newtheorem{notation}[theorem]{Notation}
\newtheorem{remark}[theorem]{Remark}
\newtheorem{observation}[theorem]{Observation}
\newtheorem{definition}[theorem]{Definition}
\newtheorem{condition}[theorem]{Condition}
\newtheorem{example}[theorem]{Example}%[section]
%\newtheoremrm{rem}{Remark}
\let\pf\proof
\let\epf\endproof
\numberwithin{equation}{section} %or {subsection}
\newcommand{\bfem}[1]{\textbf{\emph{#1}}}

\maketitle

\begin{abstract}
 In this paper, we formally define the concept of \emph{$k$-shifted contact structures} on  derived (Artin) stacks and study their local properties in the context of  derived algebraic geometry. In this regard, for  negatively shifted  contact derived $\mathbb{K}$-schemes, we   develop 
 a \emph{Darboux-like theorem} and formulate the notion of \emph{symplectification}. 
\end{abstract}
\tableofcontents

    \section{Introduction and summary}
    
    Derived algebraic geometry (DAG) essentially provides a new setup to deal with non-generic situations in geometry (e.g. non-transversal intersections and ``bad" quotients). To this end, it combines higher categorical objects and homotopy theory with many tools from homological algebra. Hence, roughly speaking, it can be considered as \emph{a higher categorical/homotopy theoretical refinement of classical algebraic geometry}. In that respect, it offers a new way of organizing information for various purposes. Therefore, it has many interactions with other mathematical domains. For a survey of some directions, we refer to \cite{Anel,Toen}.

    In the context of DAG, it is also possible to work with familiar geometric structures, but in more general forms. For instance, $k$-shifted versions of Symplectic and Poisson geometries have already been described and studied in \cite{PTVV,CPTTV}. In this regard, \cite{Brav,JS,BenBassat etall} offer some applications and local constructions.

    Throughout this paper, we mainly work within the context of Toën \& Vezzosi's version of DAG \cite{Toen,ToenHAG}. We also benefit from Lurie's version \cite{Luriethesis}. In that respect, we always consider objects with higher structures in a functorial perspective, and we focus on nice representatives for those structures. For instance, by a \emph{derived $\K$-stack}, we essentially mean a simplicial presheaf on the  category of commutative differential graded $\mathbb{K}$-algebras (cdga) having nice local-to-global  properties. %that preserves weak equivalences and possesses the descent property. 
   % \newpage
    
    DAG  provides an appropriate concept of a \emph{spectrum functor} $\spec$ from cdgas to (higher) spaces.
    Using this functor, we call a derived space of the form ${\bf X} \simeq \spec A$ for some cdga $A$ an \emph{affine derived $\K$-scheme.} As in the classical theory, a general \textit{derived $\K$-scheme} $\bf Y$ is defined to be a space which is locally modeled on ${\bf X} \simeq \spec A$. Note that  affine derived schemes are in fact the main objects of interest for us because the concepts to be discussed in this paper are all about the \emph{local structure} of derived schemes. Thus, it is enough to consider the affine case. More details  will be given in $\S$ \ref{section_prelim}.
    
    Regarding certain geometric structures on higher spaces; such as $k$-shifted (closed) $p$-forms in the sense of \cite{PTVV}, it is also known that for sufficiently ``nice" cdgas (to be clear later), we can use the $A$-module $\Omega_A^1$ of K\"{a}hler differentials as a model for the cotangent complex $\mathbb{L}_A$ of $A$ so that we write $\mathbb{L}_A\simeq \Omega_A^1$. Then, by \textit{a $k$-shifted $p$-form on $\spec A$}  for $A$ a sufficiently nice cdga, we actually mean a $k$-cohomology class of the complex $(\Lambda^p\Omega_A^1, d)$. Likewise,  a \textit{$k$-shifted closed $p$-form on $\spec A$} is just a $k$-cohomology class of the complex $\prod_{i\geq0} \big(\Lambda^{p+i} \Omega^1_{A}[-i],\ d_{tot}=d+d_{dR}\big)$. 
     
     A reasonable notion of \textit{non-degeneracy} is also available in this  framework, which leads to the definition of a \textit{shifted symplectic structure}%(cf. $\S$ \ref{secion_PTVV's symplectic geometry on spec A})
     . Loosely speaking, we are then able to define the notion of a \emph{$k$-shifted contact form} on $\spec A$ to be a $k$-shifted 1-form $\alpha$ on $\spec A$   with the property that the $k$-shifted 2-form $d_{dR}\alpha$ satisfies a non-degeneracy condition, which will be formulated later. % in $\S$ \ref{section_shifted contact structures and Darboux forms}.
      In fact, we  will provide the \emph{general definition of  contact data for derived (Artin) stacks}; rather than just for affine  derived $\K$-schemes with ``nice" local models.

     For shifted symplectic structures on derived schemes, it has been  shown in \cite{Brav} that for $k<0$, every $k$-shifted symplectic derived $\K$-scheme $(\bf X, \omega')$ is Zariski locally equivalent to $(\spec A, \omega)$ with a pair $A,\omega$  in certain symplectic Darboux form. More precisely, for $k<0$, Bussi, Brav and Joyce \cite[Theorem 5.18]{Brav} proved that given a $k$-shifted symplectic derived $\K$-scheme $(\bf X, \omega')$, one can find the so-called ``minimal standard form" cdga $A$, a Zariski open inclusion $\iota: \spec A \hookrightarrow \bf X$, and $ `` $coordinates" $x^{-i}_j, y^{k+i}_j \in A$ with $\iota^*(\omega') \sim (\omega^0, 0, 0, \dots)$ such that  \begin{equation*}
     \omega^0= \sum_{i,j} d_{dR}x_j^{-i} d_{dR}y_j^{k+i}.
     \end{equation*}We should point out that the expression of $ \omega^0 $ holds true  for the case where $k<0$ is an \emph{odd} integer. The  other possible cases require some modifications depending on whether $k/2$ is even or odd. However, the underlying idea behind the proofs for each case is the same.

     Note also that the case $k<0$ odd is relatively simple and instructive enough to capture the essential techniques for the constructions of local models under consideration. Therefore, in this paper, we will mainly concentrate on the case with $k<0$ odd and use it as a \textit{prototype construction}. For the other cases, we will not give all the details. Instead, we will only provide a brief outline. For details, we will always refer to  \cite[Examples 5.8, 5.9 \& 5.10]{Brav}. 
          
\paragraph{Results and the outline.} In this paper, we  introduce the notion of \textit{a $k$-shifted contact structure} on a derived (Artin) stack. The goals are then to develop \emph{Darboux-type models} for negatively shifted contact structures and investigate further possible outcomes, such as \emph{symplectifications}.  

A \emph{brief definition.} A $k$-shifted contact data  will consist of a  morphism  $f:\mathcal{K}\rightarrow \mathbb{T}_{{\bf X}}$ of perfect complexes, a line bundle $L$ such that $Cone(f)\simeq L[k]$, and a locally defined $k$-shifted 1-form $\alpha: \mathbb{T}_{{\bf X}} \rightarrow \mathcal{O}_{{\bf X}}[k]$ satisfying a non-degeneracy condition (cf. Definitions \ref{defn_preshiftedcontact} \& \ref{defn_shiftedcontact}). The next two theorems summarize the main results of this paper.

A \emph{Darboux-type theorem}. We first discuss the existence of Darboux-type local models for $k$-shifted contact derived $\K$-schemes when $k<0$ (see Example \ref{model example}). Next, for a locally finitely presented derived $\K$-scheme $\bf X$ with a $k$-shifted contact structure for $k<0,$ we prove the following result (cf. Theorem \ref{contact darboux}):    
    \begin{theorem} \label{THM1}
    	Every $k$-shifted contact derived $\K$-scheme $\bf X$ is locally equivalent to $(\spec A, \alpha_0)$ for $A$ a minimal standard form cdga and $\alpha_0$  a contact Darboux form.
    \end{theorem}
The \emph{symplectification}. Secondly, we establish a  shifted version of the classical connection between contact and symplectic geometries. In classical contact geometry, for a contact manifold $M$, there is a unique \emph{symplectified space} with a symplectic structure canonically determined  by  the contact data of $M$. %a standard result: One can construct a symplectified space $\widetilde{M}$ which can be seen as a certain fiber bundle over $M$ such that it equips with a unique symplectic  determined  by the contact structure of the underlying manifold $M$.   We then call $\widetilde{M}$ the \emph{symplectization} of $M.$
In this paper, we provide a similar result for shifted contact derived schemes. 

The upshot is that given a (locally finitely presented) derived $\K$-scheme $\bf X$ with a $k$-shifted contact structure for $k<0,$  we define its \emph{symplectification} $\mathcal{S}_{\bf{X}}$ of $\bf X$ as the total space  of a certain $\mathbb{G}_m$-bundle over $\bf X$, constructed via the data of $k$-shifted contact structure, and provided with a canonical $k$-shifted symplectic structure on $\mathcal{S}_{\bf{X}}$ for which the $\mathbb{G}_m$-action is of weight 1 (cf. Definition \ref{defn_symplectification} \& Theorem \ref{thm_Symplectization}). In brief, we have: %write $\mathcal{S}_{\bf{X}}(k)$ for the space of all $k$-shifted defining contact forms on $ \bf X.$ We then obtain the following theorem/definition (cf. Theorem \ref{thm_Symplectization}).
\begin{theorem} \label{THM2}
	The space $\mathcal{S}_{\bf{X}}$ has the structure of a $k$-shifted symplectic derived stack with a symplectic form $  \omega $ which is canonically determined by the shifted contact structure of $\bf X$. We then call the pair $(\mathcal{S}_{\bf{X}}, \omega)$  the \emph{symplectification of $ \bf X $}.
\end{theorem}
%\newpage
    Now, let us describe the content of this paper in more detail and provide an outline.  In Section \ref{section_review of shifted symplectic strcs}, we review derived symplectic geometry and symplectic Darboux forms.  %Section \ref{section_prelim} provides some key ideas from  To\"{e}n and Vezzosi's theory of derived algebraic geometric  and PTVV's shifted symplectic geometry: 
    We begin by some background material on derived algebraic geometry and present nice local models for derived $\K$-schemes and their cotangent complexes. In $ \S $ \ref{secion_PTVV's symplectic geometry on spec A}, using these nice local models, we %introduce (closed) $p$-forms of degree $k$ and 
    study shifted symplectic structures.
    $ \S $ \ref{the pair} outlines symplectic Darboux forms on derived schemes and presents Darboux-type results  given by Bussi, Brav and Joyce \cite[Theorem 5.18]{Brav}.

    Section 3 discusses contact structures. In $\S$ \ref{section_classical contact geo},  classical contact geometry is briefly revisited, and then, in $\S$ \ref{section_shifted contact structures and Darboux forms}, we introduce  \textit{shifted contact structures} and discuss their properties. In $\S$ \ref{section_darboux theorem}, we state a Darboux-type theorem for negatively shifted contact structures on derived $\K$-schemes  (Theorem \ref{contact darboux}) and provide the proof of Theorem \ref{THM1}.  
    
    In Section \ref{section_symplectization}, we discuss the concept of  \textit{symplectification} for  (negatively) shifted contact derived $\K$-schemes and give the proof of Theorem \ref{THM2} (cf. Definition \ref{defn_symplectification} \& Theorem \ref{thm_Symplectization}). 
    
    Section \ref{section_Artin stacs} provides some concluding remarks on possible ``stacky" generalizations of the main results of this paper. It also advertises possible future directions and some ongoing projects.

     \paragraph {Acknowledgments.} I would like to thank Alberto Cattaneo and Ödül Tetik for useful discussions.  It is also a pleasure to thank the Institute of Mathematics,  University of Zurich, where this research was conducted. I personally benefited a lot from hospitality and research environment of the Institute. 
     
     I thank the anonymous Referees for the comprehensive review. I would like to express my gratitude to the Referees for the valuable comments and suggestions, which improved the quality of the manuscript.
     
     The author acknowledges support of the Scientific and Technological Research Council of Turkey (T\"{U}BİTAK) under 2219-International Postdoctoral Research Fellowship Program (2021-1). 
    
     \paragraph{Conventions.} Throughout the paper, $ \mathbb{K} $ will be an algebraically closed field of characteristic zero. All cdgas will be graded in nonpositive degrees and  over $ \mathbb{K}.$ All classical $ \K $-schemes will be \emph{locally of finite type}, and all derived $ \K $-schemes/Artin stacks $ \textbf{X} $ are assumed to be \emph{locally finitely presented.}

    \section{Shifted symplectic structures} \label{section_review of shifted symplectic strcs}
    \subsection{Some derived algebraic geometry} \label{section_prelim}
    In this section, we outline the basics of  DAG, present some material relevant to this paper, and state some useful results from shifted symplectic geometry. As stressed before, we use both Toën \& Vezzosi's version of DAG \cite{Toen,ToenHAG} and the Lurie's version \cite{Luriethesis}. In what follows, %we will not give all the details on the theory of DAG. Instead,  
    we just  intend to give a brief sketch for the objects and constructions that we will be mostly interested in.
    \begin{definition}
    	Denote by $ cdga_{\K}$  the category of commutative differential graded $\mathbb{K}$-algebras in non-positive degrees, where an \textit{object} $A$ in $ cdga_{\K}$ consists of \begin{enumerate}
    		\item [(1)] a collection of $\K$-vector spaces $ \{A^i\}$, where $A^i$ is the $\K$-vector space of degree $i$ elements for $i=0, -1, \dots$, 
    		\item [(2)] a $\K$-bilinear, associative, supercommutative multiplication $A^n \otimes A^m \xrightarrow{\cdot} A^{n+m}$, and
    		\item  [(3)] a unique square-zero derivation of degree 1  (the differential) $\rm d$  on $A$ satisfying the graded Leibniz rule
    	\begin{equation*}
    		\rm d (a \cdot b)= (\rm d a) \cdot b + (-1)^{n}a \cdot (\rm d b) 
    	\end{equation*} for all $a\in A^n, b\in A^m.$
    	\end{enumerate} We denote such objects by $(A, \rm d)$ or just $A$. Note that $A$ has a decomposition $A=\bigoplus_i A^i.$ 
    
    A \textit{morphism} in $ cdga_{\K}$, on the other hand, is a collection of degree-wise $\K$-linear morphisms $f=\{f^i\}:A \rightarrow B$ such that each $f^i:A^i \rightarrow B^i$  commutes with all the structures of $A,B$. 
    \end{definition}

\begin{definition}
	Denote by $dSt_{\mathbb{K}}$ the $\infty$-category of derived stacks, where  an object $ \textbf{X} $ of $dSt_{\mathbb{K}}$ is given as a certain $\infty$-functor
	\begin{equation}
		\textbf{X}: \ cdga_{\mathbb{K}} \longrightarrow \ sSets,
	\end{equation}
	where $sSets$  denote the $\infty$-category of simplical sets. More precisely, objects in $dSt_{\mathbb{K}}$ are \emph{simplicial presheaves preserving weak equivalences and possessing the descent/local-to-global property w.r.t. the site structure on the source}. For a brief review, we refer to \cite{Vezz2}. 
\end{definition}
We write $cdga_{\K}^{\infty}$ for the associated $\infty$-category of 
$cdga_{\K}$ such that the homotopy category $Ho(cdga_{\K}^{\infty})$ can be obtained from $cdga_{\K}$ by formally inverting quasi-isomorphsims. 

Note that $Ho(cdga_{\K}^{\infty})$ is just an ordinary category. We should also point out that $cdga_{\K}^{\infty}$, $cdga_{\K}$ and  $Ho(cdga_{\K}^{\infty})$ have the same objects; however, lifting properties of morphisms are different. That is, a morphism $f:A\rightarrow B$ in $cdga_{\K}$ is also a morphism in $cdga_{\K}^{\infty}$ and  $Ho(cdga_{\K}^{\infty})$. But in general, the converse is not true unless $A$ is cofibrant. In the rest of this paper, we will be interested in certain types of cdgas, called \emph{standard form cdgas,} which are in fact ``sufficiently cofibrant", and hence suitable for our purposes.

In this framework, there also exists an appropriate concept of a \emph{spectrum functor} \cite[$ \S $ 4.3]{Luriethesis}
\begin{equation*}
	\spec: (cdga_{\K}^{\infty})^{op} \rightarrow dSt_{\K},
\end{equation*} which leads to the following definition.

    \begin{definition}
    	An object $\textbf{X}$ in $dSt_{\mathbb{K}}$ is called an \emph{affine derived $\mathbb{K}$-scheme} if  $\textbf{X}\simeq \spec A $ for some cdga $A \in cdga_{\K}$. An object $\bf X$ in $dSt_{\mathbb{K}}$ is then called a \emph{derived $\mathbb{K}$-scheme} if it can be covered by Zariski open affine derived $\mathbb{K}$-schemes $Y \subset X$.
    \end{definition}
    Denote by $dSch_{\K} \subset dSt_{\mathbb{K}}$ the full $\infty$-subcategory of  derived $\mathbb{K}$-schemes, and we simply write $dAff_{\K} \subset dSch_{\mathbb{K}}$ for the full $\infty$-subcategory of  affine derived $\mathbb{K}$-schemes. Note that $ \spec : (cdga_{\K}^{\infty})^{op} \rightarrow dAff_{\K} $ gives an equivalence of $\infty$-categories.

    We should note that throughout this paper, $ \mathbb{K} $ will be an algebraically closed field of characteristic zero. We also assume that all classical $ \K $-schemes are \emph{locally of finite type}, and all derived $ \K $-schemes $ \textbf{X} $ are  \emph{locally finitely presented}, by which we mean that $ \textbf{X} $ can be covered by Zariski open affines $\spec A,$ where $A$ is a finetely presented cdga over $\K.$
    
     \begin{remark} \label{rmks on Yoneda's lemma}
    	Thanks to the Yoneda embedding, one can also realize algebro-geometric objects (like classical $\K$-schemes, stacks, derived "spaces", etc...) as \textit{certain functors} in addition to the standard ringed-space formulation. We have the following enlightening diagram from \cite{Vezz2} encoding such a functorial interpretation: 
    	\begin{center}
    		\begin{tikzpicture}
    			\matrix (m) [matrix of math nodes,row sep=2.5em,column sep=7em,minimum width=2 em] {
    				 CAlg_{\mathbb{K}}   & Sets  \\
    				&  Grpds \\
    				cdga_{\mathbb{K}} &  Ssets \\};
    			\path[-stealth]
    			(m-1-1) edge  node [left] { } (m-3-1)
    			edge  node [above] {{\small schemes}} (m-1-2)
    			(m-1-1) edge  node [below] {} node [below] {{\small stacks}} (m-2-2)
    			(m-1-1) edge  node [below] {} node [below] {{\small  higher stacks}} (m-3-2)
    			(m-3-1) edge  node  [below] {{\small derived stacks}} (m-3-2)
    			
    			(m-1-2) edge  node [right] { } (m-2-2)
    			(m-2-2) edge  node [right] { } (m-3-2);
    			%edge [dashed,-] (m-2-1);
    		\end{tikzpicture}
    	\end{center} Here $ CAlg_{\mathbb{K}} $ denotes the category of commutative $\K$-algebras. Denote by $St_{\K}$ the $\infty$-category of (higher) $\K$-stacks, where objects in $St_{\K}$ are defined via the diagram above. 
    
    In the underived setup, we  have the classical ``spectrum functor" 
    \begin{equation*}
    	\cspec : (CAlg_{\mathbb{K}})^{op} \rightarrow St_{\K}.
    \end{equation*} We then call an object $X$ of $ St_{\K} $ an \emph{affine $\K$-scheme} if $X\simeq \cspec A$ for some $A\in CAlg_{\mathbb{K}}$, and a \emph{$\K$-scheme} if it has an open cover by affine $\K$-schemes.
    \end{remark}
    
    In addition to the spectrum functors $\spec, \cspec$ above, there is a natural \emph{truncation functor} $\tau: dSt_{\K} \rightarrow St_{\K}$, along with a fully faithfull left adjoint \emph{inclusion functor} $\iota: St_{\K} \hookrightarrow dSt_{\K}$, which can be thought of as an embedding of classical algebraic $\K$-spaces into derived spaces. 
    
    Note that, for a cdga $A$ there exists an equivalence $\tau \circ \spec A \simeq \cspec H^0(A)$. This means that if $\bf X$ is a (affine) derived $\K$-scheme, then its truncation $X=\tau(\textbf{X})$ is a (affine)  $\K$-scheme. Therefore, we can consider a derived $\K$-scheme $\bf X$ as \emph{an infinitesimal thickening of its truncation} $X$. It  follows that points of a derived $\K$-scheme $\bf X$ are the same as points of  of its truncation $X$. It means that the main difference between $X$ and $\bf X$ is in fact encoded by the scheme structure, not by the points!

\paragraph{Nice local models for derived $\K$-schemes.}The following result (Theorem \ref{localmodelthm}) plays an important role in constructing useful local algebraic models for  derived $\K$-schemes and for the so-called \textit{$ k $-shifted symplectic structures} on them. 

The upshot is that given a derived $\K$-scheme $\bf X$ (locally of finite presentation) and a point $x\in \bf X$, one can always find a ``refined" local affine neighborhood $ \spec A $ of $x$ that allows us to make more explicit computations over this neighborhood. For example, using such local models, we can identify the cotangent complex $\mathbb{L}_A$ with the module of K\"{a}hler differentials $\Omega_A^1$, and then we can provide explicit representatives (rather than just cohomology classes) for (closed) $p$-forms of degree $k$.  In this regard, Bussi, Brav and Joyce proved the following theorem.
    
    \begin{theorem} (\cite[Theorem 4.1]{Brav})\label{localmodelthm}
    	Every derived $\mathbb{K}$-scheme $X$ is Zariski locally modelled on $ \spec A $ for some ``minimal standard form" cdga $ A $ in $cdga_{\mathbb{K}}$.
    \end{theorem}
 More precisely,  for each $x\in X$ there is a pair $\big(A, i: \spec A \hookrightarrow X\big)$  and $p \in \spec H^0(A)$ such that $i$ is an open inclusion with $i(p)=x$, where $A$ is a special kind of cdga (cf. Definition \ref{defn_standard forms}). 
 
 Moreover, there is a reasonable way to compare two such local charts $i: \spec A \hookrightarrow X$ and $j: \spec B \hookrightarrow X$ on their overlaps via a third chart. For details, see \cite[Theorem 4.1 \& 4.2]{Brav}.

In the remainder of this section, we shall elaborate the content of Theorem \ref{localmodelthm}, and introduce appropriate notions for the constructions of interest. We will closely follow \cite{Brav,JS}.
    
\begin{definition} \label{defn_standard forms}
	$ A \in cdga_{\mathbb{K}}$ is of \emph{standard form} if  $A^0$ is a smooth finitely generated $\mathbb{K}$-algebra, the  module $\Omega^1_{A^0}$ of K\"{a}hler differentials is free $A^0$-module of finite rank, and the graded algebra $A$ is freely generated over $A^0$ by finitely many generators, all in negative degrees.

\end{definition}    
  
  In fact, there is a systematic way of constructing such cdgas. \cite[Example 2.8]{Brav} explains how to build these cdgas starting from a smooth $\K$-algebra $A^0:=A(0)$ via applying a sequence of localizations. The upshot is follows: Let $n\in \N$, then a cdga $A$, as a commutative graded algebra, can be constructed inductively from a smooth $\K$-algebra $A(0)$ % (with $\Omega^1_{A^0}$  free $A^0$-module of finite rank) 
  by adjoining free finite rank modules $M^{-i}$ of generators in degree $-i$ for $i=1,2,\dots, n.$ More precisely, for any given $n\in \N$, we can inductively construct a sequence of cdgas 
  \begin{equation} \label{A(n) construction}
  	A(0) \rightarrow A(1) \rightarrow \cdots \rightarrow A(i)\rightarrow \cdots A(n)=:A,
  \end{equation} where  $ A^0:=A(0) $, and $A(i)$ is obtained from $A(i-1)$ by adjoining generators in degree $-i$, given by $M^{-i}$, for all $i$. Here, each $M^{-i}$ is a free finite rank module (of degree $-i$ generators) over $A(i-1)$. Therefore, the underlying commutative graded algebra of $A=A(n)$ is freely generated over $A(0)$ by finitely many generators, all in negative degrees $-1,-2,\dots, -n$.

 \begin{definition} \label{1st defn of minimalty}
 	A standard form cdga $ A $ is said to be \emph{minimal} at  $p \in \spec H^0(A)$ if $ A=A(n) $ is defined by using the minimal possible  numbers of graded generators in each degree $\leq 0$ compared to all other cdgas locally equivalent to  $ A $ near $p$. (There will be an equivalent definition for minimality later, see Definition \ref{2nd defn of minimalty}.)
 \end{definition}  

 \begin{definition}
 	Let $A$ be a standard form cdga.  $A'\in cdga_{\mathbb{K}}$ is called a \emph{localization} of $A$ if $A'$ is obtained from $A$ by inverting an element $f\in A^0$, by which we mean $A'=A\otimes_{A^0}A^0[f^{-1}]$. %for $f\in A^0$. 
 	
 	$A'$ is then of standard form with $A'^0 \simeq A^0[f]$. If $p\in \cspec H^0(A)$ with $ f(p)\neq0 $, we say $A'$ is a \emph{localization of $A$ at $p$}.
 \end{definition}
 With these definitions in hand, one has the following observations: 
 \begin{observation}
		Let $A$ be a standard form cdga. If $A'$ is  a \emph{localization} of $A$, then $\spec A' \subset \spec A$ is a Zariski open subset. Likewise, if $A'$ is a \emph{localization of $A$ at $p\in \cspec H^0(A)\simeq \tau(\spec A)$}, then  $\spec A' \subset \spec A$ is a Zariski open neighborhood of $p.$
	\end{observation}
\begin{observation}
	 Let  $ A=A(k)$ be a standard form cdga, then there exist generators $x_1^{-i}, x_2^{-i}, \cdots, x_{m_i}^{-i} $ in $A^{-i}$ (after localization, if necessary) with $i= 1, 2, \cdots, k$  \ and $m_i \in \mathbb{Z}_{\geq 0}$ such that 
	\begin{equation}
	A = A(0) \big[x_j^{-i} : i= 1, 2, \dots, k, \ j= 1,2, \dots, m_i\big],
	\end{equation} where the subscript $j$ in $x_j^i$ labels the generators, and the superscript $i$ indicates the degree of the corresponding element.  So, we can consider  $A$ as a \emph{graded polynomial algebra over $A(0)$ on finitely many generators, all in negative degrees.} 
\end{observation}%Moreover, since $A^0$ is finitely generated $\mathbb{K}$-algebra, there exists a surjection $\mathbb{K}[x_1^0, \cdots, x_{m_0}^0] \twoheadrightarrow A^0$ such that $A^0= \mathbb{K}[x_1^0,x_2^0, \cdots, x_{m_0}^0]/ I$ for some ideal $I$. %We shall be interested in the cases where the ideals $I$ are also finitely generated. Then $A^0$ is indeed said to be \emph{finitely presented}.
%	\item In general, we are interested in minimal standard form cdgas for which each number $m_i$ above is taken to be the least possible number such that $m_i=dim (H^i (\mathbb{L}_{A} |_p))$  and $d^i |_p =0$ for $i= -1, -2, \cdots, k$, and $p \in spec (H^0(A))$. 
\begin{definition}
	We then define the \emph{virtual dimension} of $A$ to be the integer $\vdim A= \sum _i (-1)^im_i$. 
\end{definition}
	\begin{observation}
		Geometrically, the ``smoothness$ " $ condition on $A^0$ implies that the corresponding affine $\mathbb{K}$-scheme $U=\cspec A^0$ is smooth together with a local ($\acute{e}tale$) coordinate system\begin{equation}
	(x_1^0,x_2^0, \dots, x_{m_0}^0): U \longrightarrow \mathbb{A}^{m_0}_{\mathbb{K}}.
	\end{equation} 
	\end{observation}

\paragraph{Nice local models for cotangent complexes of derived schemes.} Given	$ A \in cdga_{\mathbb{K}}$, $ d$ on $A$ induces a differential on $\Omega_A^1$, denoted again by $ d$. This makes $\Omega_A^1$ into a $\rm dg$-module $(\Omega_A^1, d)$ with the property that $\delta \circ  d = d \circ \delta,$ where $\delta: A \rightarrow \Omega_A^1$ is the universal derivation of degree $0$. 

Write  the decomposition of $ \Omega^1_{A} $ into graded pieces 
$ \Omega^1_{A} = \bigoplus_{k=-\infty}^0 \big(\Omega^1_{A}\big)^k $
with the  differential $d: \big(\Omega^1_{A}\big)^k \longrightarrow \big(\Omega^1_{A}\big)^{k+1}$. Then we define the \textit{de Rham algebra of $A$} as a  double complex 
   \begin{equation}
   DR(A)= Sym_A(\Omega_A^1[1]) \simeq \bigoplus \limits_{p=0}^{\infty} \bigoplus \limits_{k=-\infty}^0 \big(\Lambda^p \Omega^1_{A}\big)^k [p],
   \end{equation} where $ DR(A) $ has two gradings: the grading w.r.t. $ p $ is called the \emph{weight}, and the grading w.r.t. $ k $ is called the \emph{degree}. By construction, there are two differentials, namely the \emph{internal differential d} and the \emph{de Rham differential} $d_{dR}.$ We diagrammatically have 
   
 %  \begin{align}
 %  d &: \big(\Lambda^p \Omega^1_{A}\big)^k [p] \longrightarrow \big(\Lambda^p \Omega^1_{A}\big)^{k+1} [p], \\ %\text{\ \ \ \ \ \ \ (\emph{increasing the degree by 1})}\\
 %  d_{dR} &: \big(\Lambda^p \Omega^1_{A}\big)^k [p] \longrightarrow \big(\Lambda^{p+1} \Omega^1_{A}\big)^k [p+1] %\text{\ \ \ \ \  \ (\emph{increasing the weight by 1})}
 %  \end{align} 

   \begin{equation} \label{deRham cmpx}
   \begin{tikzpicture}
   \matrix (m) [matrix of math nodes,row sep=1em,column sep=2em,minimum width=1em] 
   { % 4x4 matrix
   	& \vdots         & \vdots &      \\
   	\cdots  & \big(\Lambda^{p+1} \Omega^1_{A}\big)^k [p+1] &\big(\Lambda^{p+1} \Omega^1_{A}\big)^{k+1} [p+1]       &\cdots\\
   	\cdots  & \big(\Lambda^p \Omega^1_{A}\big)^k [p] & \big(\Lambda^p \Omega^1_{A}\big)^{k+1} [p]   &\cdots \\
   	& \vdots         & \vdots &        \\
   };
   \path[-stealth]
   (m-2-2) edge  node [left] {{\small $  $}} (m-1-2)
   (m-3-2) edge  node [left] {{\small $ d_{dR} $}} (m-2-2)
   (m-4-2) edge  node [left] {{\small $  $}} (m-3-2)
   %edge  node [above] {{\small $ \Phi_V $}} (m-1-2)
   (m-2-1.east|-m-2-2) edge  node [below] {} node [below] {{\small $  $}} (m-2-2)
   (m-3-1.east|-m-3-2) edge  node [below] {} node [below] {{\small $  $}} (m-3-2)
   
   (m-2-2.east|-m-2-3) edge  node [below] {} node [above] {{\small $ d $}} (m-2-3)
   (m-3-2.east|-m-3-3) edge  node [below] {} node [below] {{\small $ d $}} (m-3-3)
   
   (m-2-3.east|-m-2-4) edge  node [below] {} node [below] {{\small $  $}} (m-2-4)
   (m-3-3.east|-m-3-4) edge  node [below] {} node [below] {{\small $  $}} (m-3-4)
   
   (m-2-3) edge  node [right] {{\small $ $}} (m-1-3)
   (m-3-3) edge  node [right] {{\small $ d_{dR}$}} (m-2-3)
   (m-4-3) edge  node [right] {{\small $ $}} (m-3-3);
   
   %edge [dashed,-] (m-2-1);
   \end{tikzpicture} 
   \end{equation} such that $d_{tot}= d + d_{dR}$, and both differentials satisfy the  relations
   \begin{equation}
   d^2=d_{dR}^2=0, \text{\ and \ } d\circ d_{dR} + d_{dR} \circ d =0.
   \end{equation} We also have the natural multiplication on $ DR(A)$:
   \begin{equation}
   \big(\Lambda^p \Omega^1_{A}\big)^k [p] \times \big(\Lambda^q \Omega^1_{A}\big)^{\ell} [q] \longrightarrow \big(\Lambda^{p+q} \Omega^1_{A}\big)^{k+\ell} [p+q].
   \end{equation} 
   
  \begin{observation}
  	 It should be noted that the constructions of $\Omega_A^1 \text{ and } DR(A)$ depend only on the underlying commutative graded algebra of $A$, not on the differential $d$ on $A.$ 
  \end{observation}

   \begin{remark} When $ A=A(k) $ is a minimal standard form cdga, there are two important outcomes: %\cite{Brav,JS}:
  	\begin{enumerate}
  		%\item []
  		 \item With such local coordinates $(x_1^0,x_2^0, \cdots, x_{m_0}^0)$, we have \begin{equation}
  		\Omega^1_{A^0} \cong A^0 \otimes_{\mathbb{K}} \langle d_{dR}x_1^0, \cdots, d_{dR}x_{m_0}^0\rangle_{\mathbb{K}}.
  		\end{equation}Furthermore, the K\"{a}hler differentials is a $A$-module of the form
  		\begin{equation}
  		\Omega^1_{A} \cong A \otimes_{\mathbb{K}}  \langle d_{dR}x_j^{-i} : i= 0,1, 2, \cdots, k, \ j= 1,2, \cdots, m_i \rangle_{\mathbb{K}}.
  		\end{equation}
  		\item $\Omega^1_{A}$ provides a local model for the cotangent complex $\mathbb{L}_{A}$. That is, in the case of a minimal standard form cdga, the \textit{cotangent complex} $\mathbb{L}_{A}$ has the identification \begin{equation} \label{cotangent cpx vs kahler diff}
  		\mathbb{L}_{A}= \Omega^1_{A}.
  		\end{equation}  
  	\end{enumerate}
  \end{remark}  Note that if $D(Mod_{A})$ denotes the derived category of $Mod_{A}$, then $\mathbb{L}_{A} \in D(Mod_{A})$ for standard form cdgas. In general, even if both $\mathbb{L}_{A} \text{ and } \Omega^1_{A}$ are closely related, the identification in (\ref{cotangent cpx vs kahler diff}) is not true for an arbitrary   $ A \in cdga_{\mathbb{K}}$ \cite{Brav}. 

When $A=A(n)$ is a standard form cdga as in (\ref{A(n) construction}), we also have the following description for the restriction of the cotangent complex $\mathbb{L}_A $ to $\cspec H^0(A)$.

 \begin{proposition} (\cite[Prop. 2.12]{Brav}) \label{proposition_L as a complex of H^0 modules}
If $A=A(n)$, with $n\in \N$, is a standard form cdga constructed inductively as in (\ref{A(n) construction}), then the restriction of $\mathbb{L}_A $ to $\cspec H^0(A)$ is represented by a complex of $H^0(A)$-modules
\begin{equation} \label{complex of H^0 modules}
0\rightarrow V^{-n} \xrightarrow {d^{-n}} V^{-n+1} \rightarrow \cdots \rightarrow V^{-1} \xrightarrow {d^{-1}} V^0 \rightarrow 0, 
\end{equation}where each $V^{-i}$ can in fact be defined as $V^{-i}= H^{-i} (\mathbb{L}_{A(i)/A(i-1)})$, with $\mathbb{L}_{A(i)/A(i-1)}$  the \emph{relative cotangent complex} of the map $A(i-1) \rightarrow A(i)$ in (\ref{A(n) construction}) satisfying \begin{equation*}
\mathbb{L}_{A(i)/A(i-1)}\simeq A(i) \otimes_{A(i-1)} M^{-i} [i].
\end{equation*}  Moreover, the differential $V^{-i} \xrightarrow {d^{-i}} V^{-i+1}$ is identified with the composition \begin{equation*}
H^{-i} (\mathbb{L}_{A(i)/A(i-1)}) \rightarrow H^{-i+1} ( \mathbb{L}_{A(i-1)}) \rightarrow H^{-i+1} (\mathbb{L}_{A(i-1)/A(i-2)}),
\end{equation*} which can be obtained from the fiber sequences induced by the maps $A(i-1) \rightarrow A(i)$ in (\ref{A(n) construction}). Note in particular that $j>-i,$ we have $ H^{j} (\mathbb{L}_{A(i)/A(i-1)})=0. $ More details and the proof can be found in \cite[Prop. 2.12]{Brav}.
 \end{proposition}   With this result in hand,  using local coordinates above, write $$V^{-i}=\langle d_{dR} x^{-i}_1, d_{dR} x^{-i}_2, \dots, d_{dR} x^{-i}_{m_i} \rangle_{A(0)}  \ \text{ for } i=0,1, \dots, n.$$It follows that we have a similar local description for the tangent complex $\mathbb{T}_A= (\mathbb{L}_A)^{\vee}$ of $A$ when restricted to  $\cspec H^0(A)$. Also, we have an alternative definition of minimality at a point $p\in \cspec H^0(A)$ for a cdga of the form $A=A(n)$.

\begin{definition} \label{2nd defn of minimalty}
Let $A=A(n)$, with $n\in \N$, be a standard form cdga constructed inductively as in (\ref{A(n) construction}). $A$ is said to be \emph{minimal at} $p\in \cspec H^0(A)$ if the internal differential $d^{-i}|_p = 0$ in the complex $\mathbb{L}_A|_{\cspec H^0(A)}$ given in (\ref{complex of H^0 modules}).
\end{definition}
    Note from Proposition \ref{proposition_L as a complex of H^0 modules} that at $p\in \cspec H^0(A)$, with $A$ a standard form cdga, $\Omega^1_{A} |_p$ is a complex of $\K$-vector spaces and that Definition \ref{2nd defn of minimalty} implies  $m_i=\dim (H^{-i} (\mathbb{L}_{A} |_p))$ for each $i$, and hence $ A $ is defined by using the minimum number of graded variables in each degree $\leq 0$ compared to all other cdgas locally equivalent to  $ A $ near $p$. Therefore, one  recovers Definition \ref{1st defn of minimalty}.

\subsection{PTVV's shifted symplectic geometry on derived schemes}   \label{secion_PTVV's symplectic geometry on spec A} 
  Let $\bf X$ be a locally finitely presented derived $\K$-scheme with $p\geq 0, \ k\in \Z$. Pantev et al. \cite{PTVV} define simplicial sets of $p$-forms of degree $k$  and closed $p$-forms of degree $k$  on $\bf X.$ Denote these simplicial sets by $\mathcal{A}^p(X,k)$ and $\mathcal{A}^{(p,cl)}(X,k)$, respectively. These definitions are in fact given first for affine derived $\K$-schemes. Later, both concepts are defined for a general $\bf X$ in terms of mapping stacks. A summary of key ideas can be found in \cite[$\S$ 3.4]{Brav}
  
  In our case, we consider $ \textbf{X} = \spec A$ with $A$ a standard form cdga, and hence take $ \Lambda^p\mathbb{L}_A=\Lambda^p \Omega^1_{A}$.  Therefore,  elements of $ \mathcal{A}^{p}(X,k)$  form a simplicial set such that $k$-cohomology classes of the complex $ \big(\Lambda^p \Omega^1_{A}, d\big) $  correspond to the connected components of this simplicial set. Likewise,  the connected components of $\mathcal{A}^{(p,cl)}(X,k)$ are identified with the $k$-cohomology classes of the complex $\prod_{i\geq0} \big(\Lambda^{p+i} \Omega^1_{A}[-i], d_{tot}\big)$. We want to work with explicit representatives for these  classes. 
  
  It should be noted that the results that are cited or to be proven in this paper are all about the \emph{local structure} of derived schemes. Thus, it is enough to consider the affine case. Moreover, we always assume all local models are sufficiently nice by using Theorem \ref{localmodelthm} if necessary.
 \begin{definition} \label{defn_p form of deg k}
 	Let $\textbf{X}=\spec A$ be an affine derived $\mathbb{K}$-scheme for $A$ a minimal standard form cdga. A \emph{$p$-form of degree $k$ on $\bf X$} for $p\geq 0$ and $k \leq 0$ is an element \begin{equation}
 	\omega^0 \in \big(\Lambda^p \Omega^1_{A}\big)^k \text{ with } d\omega^0=0. 
 	\end{equation}
 \end{definition} \noindent Note that an element  $\omega^0$ defines a cohomology class as being $d$-closed. That is, $$[\omega^0] \in H^k\big(\Lambda^p \Omega^1_{A}, d\big),$$ where two $ p $-forms  $\omega_1^0, \omega_2^0 $ of degrees $k$ are \emph{equivalent} if there exists $ \alpha_{1,2} \in \big(\Lambda^p \Omega^1_{A}\big)^{k-1}$ so that $\omega_1^0-\omega_2^0= d\alpha_{1,2}$.

    \begin{remark} \label{rmk_underived setup}
    	In the classical $ `` $underived$ " $ case, for instance when $X=\cspec A$ is smooth for a commutative $\mathbb{K}$-algebra $A$, the cotangent complex $\mathbb{L}_X$ is just a vector bundle over $X$, and denoted simply by $T^*X$.   Then, a $p$-form $\omega$ on $ X $ is defined to be a global section of the bundle $\Lambda^pT^*X$. A careful observation reveals that Definition \ref{defn_p form of deg k} does generalize the definition of a $p$-form on a smooth space in the following sense: It is clear that any commutative $\mathbb{K}$-algebra $A$ can be realized as an object in $cdga_\mathbb{K}$ concentrated in degree 0 with the trivial differential. Thus, in the language of Definition \ref{defn_p form of deg k}, a na\"{\i}ve notion of $ `` $a $p$-form $\omega$ on a smooth space $ X "$ is just a $p$-form $\omega$ of degree 0 on $ X\simeq\tau\circ \iota(X)$ in $St_{\K}$ such that $ \omega \in \big(\Lambda^p T^*X\big)^0.  $ 
    	Note that the condition $d\omega=0$ holds trivially, and hence $ [\omega] \in H^0\big(\Lambda^p T^*X, d=0\big)$. Here, $\mathbb{L}_X=T^*X$ is again viewed as graded object concentrated in degree 0, with the zero differential.
  \vspace{5pt}
    	
    	In DAG, on the other hand, $ \Lambda^p \mathbb{L}_X $ is a (double) complex which possesses a non-trivial internal differential as above, and hence one needs to take into account higher non-trivial cohomology groups as well.
    \end{remark}

 \begin{definition} \label{defn_shifted non-degeneracy}
 	A 2-form $\omega^0$ of degree $ k $ on $\spec A $ for $A$ a minimal standard form cdga is called \emph{non-degenerate} if the induced morphism $ \omega^0: \mathbb{T}_{A} \rightarrow \Omega^1_{A}[k], \ Y \mapsto \iota_{Y} \omega^0$, is a quasi-isomorphism, where $\mathbb{T}_{A}=(\mathbb{L}_{A})^{\vee}=\Hom_{A}(\Omega^1_{A},A)$ is \textit{the tangent complex of $A$.}
 \end{definition}   
    
   \begin{definition}\label{defn_closed p form}
   Let $\textbf{X}=\spec A $ be an affine derived $\mathbb{K}$-scheme with $A$ a minimal standard form cdga. A \emph{closed $p$-form of degree $k$ on $\bf X$} for $p\geq 0$ and $k \leq 0$ is a sequence $\omega=(\omega^0, \omega^1, \dots)$ with $\omega^i \in \big(\Lambda^{p+i} \Omega^1_{A}\big)^{k-i}$ satisfying the following conditions:
   \begin{enumerate}
   	\item [(1)] $d\omega^0=0$ in $ \big(\Lambda^p \Omega^1_{A}\big)^{k+1}.$
   	\item [(2)] $d_{dR}\omega^i + d\omega^{i+1}=0$ in $\big(\Lambda^{p+i+1} \Omega^1_{A}\big)^{k-i}$, $i \geq 0.$
   \end{enumerate}
   \end{definition}

 \begin{remark} 
 	\begin{enumerate}
 		\item [ ]
 		\item From Definition \ref{defn_closed p form}, there exists a natural projection morphism \begin{equation}
 		\pi: \mathcal{A}^{(p,cl)}(X,k)\longrightarrow \mathcal{A}^{p}(X,k), \ \ \omega=(\omega^i)_{i\geq 0} \longmapsto \omega^0.
 		\end{equation}
 		
 		\item When we restrict ourselves to the classical case as in Remark \ref{rmk_underived setup}, \emph{the one in which everything is concentrated in degree 0}, we have $d=0$ and hence $d_{tot}=d_{dR}$. Moreover, the only possible non-trivial component of $\omega$ is $\omega^0$. Therefore, using the truncation functor as before, the conditions in Definition \ref{defn_closed p form} reduce to  \begin{equation}
 		\omega^0 \in H^0\big(\Lambda^p T^*X\big) \text{  with } d_{dR}\omega^0=0.
 		\end{equation} Thus, Definition \ref{defn_closed p form} reduces to the usual definition of a (de Rham) closed $p$-form on smooth spaces.
 	\end{enumerate}
 	 
 \end{remark}

\begin{definition}
	A closed $2$-form $ \omega=(\omega^i)_{i\geq 0} $ of degree $k$ on  an affine derived $\mathbb{K}$-scheme $\spec A $ for a minimal standard form cdga $A$ is called \emph{a $k$-shifted symplectic structure} if $\pi(\omega)=\omega^0$ is a non-degenerate 2-form of degree $k$.
\end{definition}

\subsection{Shifted symplectic Darboux models} \label{the pair}
One of the main theorems in  \cite{Brav} provides a $ k $-shifted version of the classical Darboux theorem in symplectic geometry. The statement is as follows. 

\begin{theorem} (\cite[Theorem 5.18]{Brav}) \label{Symplectic darboux}
	Given a derived $\mathbb{K}$-scheme $X$ with a $k$-shifted symplectic form $\omega'$ for $k<0$ and $x\in X$, there is a local model $\big(A, f: spec A \hookrightarrow X, \omega \big)$  and $p \in \spec H^0(A)$ such that $f$ is an open inclusion with $f(p)=x$, $A$ is a standard form that is minimal at  $p$, and $\omega$ is a $k$-shifted symplectic form on $\spec A$ such that $A, \omega$ are in Darboux form, and $f^*(\omega') \sim \omega$  in the space of $k$-shifted closed 2-forms.
\end{theorem} 

To be more precise, it is proven in \cite[Theorem 5.18]{Brav} that such  $ \omega $  can be constructed explicitly depending on the integer $k<0$. Indeed, there are three cases in total: 
\begin{equation*}
(1) \ k \text{ is odd}, \ (2) \ k\equiv 0 \mod 4, \ (3) \ k\equiv 2 \mod 4 .
\end{equation*}
Equivalently, the cases can be expressed as $ (1) \ k/2\notin \Z, \ (2) \ k/2 \text{  is even}, \text{ and } (3) \ k/2 \text{  is odd},$ respectively. In short, Theorem \ref{Symplectic darboux} says that every $k$-shifted symlectic derived $\K$-scheme $(\bf X, \omega')$ is Zariski locally equivalent to $(\spec A, \omega)$ for some $A,\omega$, where  $A$ is a minimal standard form cdga and $\omega$ is a $k$-shifted symplectic form on $\spec A$ such that $\omega$ is given in a standard way depending on the cases above. 

In this paper, for simplicity, we will examine a family of explicit Darboux forms for   $k<0$ an \emph{odd} integer \cite[Example 5.8]{Brav}. The other cases can be studied in a similar way, but with some modifications.  We will outline the steps. More details can be found in \cite[$\S$ 5.3]{Brav}. %and, in particular, $k=-1$.

%Before discussing contact geometry and contact local models on $\spec A$, 
We first begin with a useful result  that plays a significant role in constructing Darboux-type local models below. The upshot is that one can always simplify the form of a given closed 2-form $ \omega= (\omega^0,\omega^1, \omega^2, \dots)$ of degree $k<0$ on $\spec A$ so that $\omega^0$ can be taken to be exact and $\omega^i=0$ for all $i>0.$ More precisely, we have the following result.

\begin{proposition} (\cite[Prop. 5.7]{Brav} )\label{Proposition_exactness}
	Let  $ \omega= (\omega^0,\omega^1, \omega^2, \dots)$ be a closed 2-form of degree $k<0$ on $\spec A$ for $A$ a standard form cdga over $\K.$ Then there exist $H\in A^{k+1}$ and $\phi \in (\Omega^1_{A})^k$ such that $dH=0$ in $A^{k+2}$, $d_{dR}H+d\phi=0$ in $(\Omega^1_{A})^{k+1}$, and $\omega \sim (d_{dR}\phi, 0, 0, \dots).$ %In fact, $d_{dR}\phi=k\omega^0.$ 
	
	Moreover, if $(H', \phi')$ is another such pair for fixed $\omega, k, A,$ then there exist $h\in A^k$ and $\sigma\in (\Omega_A^1)^{k-1}$ such that $H-H'=dh$ and $\phi - \phi'=d_{dR}h+d\sigma.$
\end{proposition}
The proof of Proposition \ref{Proposition_exactness} is based on the fact that any such forms can be interpreted in the context of \emph{cyclic homology theory of mixed complexes}. Indeed, any such forms can be viewed as cocycles in the so-called \emph{negative cyclic complex of weight p} on $\spec A$, which is constructed from the de Rham algebra $DR(A)$ in certain way. When $p=2$, there are some useful short exact sequences and vanishing results, by which one can eventually obtain the desired simplification above. For more details on this cyclic homology perspective, we  refer to \cite[$\S$ 5.2]{Brav}.

\begin{observation} \label{observation_making H exact}
	Assume $(H, \phi)$ is a such pair for fixed $\omega, k, A,$ with $d_{dR}\phi=k\omega^0.$ Let $f\in \K$ be a non-zero element. Define $H'=fH$ and $\phi'=f\phi$. Then both $H',\phi'$ satisfy the relations $dH'=0$ and $d_{dR}H'+d\phi'=0$. From the choices, we also have $d_{dR}\phi'=kf\omega^0,$ and hence  $\omega \sim (d_{dR}\phi', 0, 0, \dots).$ By Proposition \ref{Proposition_exactness}, there exist $h\in A^k$ and $\sigma\in (\Omega_A^1)^{k-1}$ such that $H-H'=dh$ and $\phi - \phi'=d_{dR}h+d\sigma.$ It follows that $(1-f)H=dh$. Localizing $A$ by the element $(1-f)$ if necessary, we can write $H=d\big[(1-f)^{-1}h\big].$ It means that we can $ ``locally" $ take $H$ to be $ d $-exact. 
\end{observation}
\paragraph{Prototype Darboux model for $k<0$ odd.}%Now, we consider the following Darboux-type local model for shifted symplectic derived $\K$-schemes as a prototype example.

 Let $k=-2\ell-1$ for $\ell \in \N$. Then the \textit{local model} consists of the following data:
\begin{enumerate}
	\item Let $A^0=A(0)$ be a smooth $\mathbb{K}$-algebra of $\dim m_0$, choose $x_1^0, \dots, x_{m_0}^0$ such that $d_{dR}x_1^0, \dots, d_{dR}x_{m_0}^0$ form a basis for $\Omega_{A^0}^1$. Then  $A$ is defined to be the free graded algebra over $A^0$ generated by variables  
	\begin{align} \label{local variables}
	 x_1^{-i}, x_2^{-i}, \dots, x_{m_i}^{-i} & \text{ in degree } (-i) \text{ for } i= 1, 2, \dots, \ell, \nonumber \\
	 y_1^{k+i}, y_2^{k+i}, \dots, y_{m_i}^{k+i} &\text{ in degree } (k+i) \text{ for } i=0,1,\dots, \ell.
	\end{align}
	\item $\Omega^1_{A}$ is  the free  $A$-module of finite rank given by
	\begin{equation}
	\Omega^1_{A} \simeq A\otimes_{\mathbb{K}}  \langle d_{dR}x_j^{-i}, d_{dR}y_j^{k+i} : i= 0, 1, \dots, \ell, \ j= 1,2, \cdots, m_i \rangle_{\mathbb{K}}.
	\end{equation}
	\item Define an element $\omega^0 \in (\Lambda^2\Omega^1_{A})^k$ of degree $k$ and weight 2 in $DR(A)$ to be \begin{equation}
	\omega^0=\displaystyle \sum_{i=0}^{\ell} \sum_{j=1}^{m_i} d_{dR}x_j^{-i} d_{dR}y_j^{k+i}.
	\end{equation}
	 \item It follows from  Proposition \ref{Proposition_exactness} that there exists a pair $(\phi, H)\in (\Omega^1_{A})^k \times A^{k+1}$ satisfying the following properties:
	 \begin{enumerate} 
	 	 
	 	\item $dH=0$ in $A^{k+2}$, \ $d_{dR}H+d\phi=0$ in $(\Omega^1_{A})^{k+1}$, and $d_{dR}\phi=k\omega^0.$ 
	 	
	 		\item $H$ satisfies the condition (a.k.a. \emph{the classical master equation})
	 	\begin{equation} \label{defn_CME}
	 	\displaystyle \sum_{i=1}^{\ell} \sum_{j=1}^{m_i} \dfrac{\partial H}{\partial x_j^{-i}} \dfrac{\partial H}{\partial y_j^{k+i}}=0 \text{ in } A^{k+2}, 
	 	\end{equation} which in fact corresponds to the condition $``dH=0".$ We call $H$ the \emph{Hamiltonian.}

 	\item Explicitly, \ we have \begin{equation} \label{defn_phi}
 	\phi:= \displaystyle \sum_{i=0}^{\ell} \sum_{j=1}^{m_i} \big[-ix_j^{-i}d_{dR}y_j^{k+i} +(k+i)y_j^{k+i}d_{dR}x_j^i\big].
 	\end{equation} Note that we can choose another representatives  by replacing $H, \phi$ by $\phi'=\phi + d_{dR}\theta $ and $H'=H + d\theta $ for any $\theta\in A^k$. This modification will leave $\omega^0$ unchanged, and both $H', \phi'$ satisfy $dH'=0$ and $d_{dR}H'+d\phi'=0$. Letting $\theta= \sum_{i=0}^{\ell} \sum_{j=1}^{m_i} \big[(-1)^ix_j^{-i}y_j^{k+i}\big]$, for instance, we may take $\phi:= k\sum_{i=0}^{\ell} \sum_{j=1}^{m_i} y_j^{k+i}d_{dR}x_j^{-i}$. 
 
 	\item The internal differential $d$ on $A$ can be defined as %$d=0$ on $A^0,$ and 
 	\begin{equation} \label{defn_internal d}
 	d|_{A^0}=0, \ \ dx_j^{-i} =  \dfrac{\partial H}{\partial y_j^{k+i}} \ \text{ and } \ dy_j^{k+i} =  \dfrac{\partial H}{\partial x_j^{-i}}.
 	\end{equation}
	 \end{enumerate}
 
	\item Clearly  $ d_{dR} \omega^0=0 $, but it is a little bit cumbersome to check that $d\omega^0=0$, and $\omega^0$ defines a non-degenerate pairing. For details, we refer to \cite[Example 5.8]{Brav}. As a result, the sequence %\begin{equation} \label{local model for omega}
	$ \omega:= (\omega^0,0,0, \dots)$ %, \text{ with } \omega^i=0 \ \ \forall i > 0, 
	%\end{equation} 
	defines a $k$-shifted symplectic structure on $\spec A$. 
\end{enumerate}
\begin{definition}
	If $A, \omega$ are as above, then we say that the pair \emph{$(A, \omega)$ is in (symplectic) Darboux form}.
\end{definition}

 In brief, Theorem \ref{Symplectic darboux} implies that every $k$-shifted symlectic derived $\K$-scheme $(\bf X, \omega')$, with $k<0$ odd, is Zariski locally equivalent to $(\spec A, \omega)$ for a pair $A,\omega$  in Darboux form as above. More precisely, Bussi, Brav and Joyce \cite[Theorem 5.18]{Brav} proved that given a $k$-shifted symlectic derived $\K$-scheme $(\bf X, \omega')$, one can find a minimal standard form cdga $A$ with $ `` $coordinates" $x^{-i}_j, y^{k+i}_j$, and a Zariski open inclusion $\iota: \spec A \hookrightarrow \bf X$ such that $\iota^*(\omega') \simeq (\omega^0, 0, 0, \dots)$ and $ \omega^0= \sum_{i,j} d_{dR}x_j^{-i} d_{dR}y_j^{k+i}.$ 
 
 Note that the expression of $ \omega^0 $ above is actually valid for the case $k\not\equiv2 \mod 4$, but the case $k\equiv2 \mod 4$ requires some modifications and extra variables. In either case, as mentioned before, the proofs follow the same logic. We now give an outline for the cases.
 \paragraph{Darboux forms for the other cases of $ k.$} For the sake of completeness, we briefly summarize the cases when $k/2$ is even or odd. Here, the main difference from the case $k$ being odd is about the existence of\textit{ middle degree variables}. In fact, when $k$ is odd, there is no such degree. But if $k/2$ is even, there are such variables and 2-forms  are \textit{anti-symmetric} in these variables. On the other hand, when $k/2$ is odd, such forms are \textit{symmetric} in the middle degree variables. Let us briefly examine each case: 
 
 \begin{enumerate}
 	\item[(a)] \cite[Example 5.9]{Brav} When $k/2$ is \emph{even}, say $k=-4\ell$ for $\ell\in \N$, the cdga $A$ is now free over $A(0)$ generated by the new set of variables 
 	\begin{align} \label{new local variables for k=-4l}
 	& x_1^{-i}, x_2^{-i}, \dots, x_{m_i}^{-i}  & & \text{ in degree } -i \text{ for } i= 1, 2, \dots, 2\ell-1, \nonumber \\
 	& x_1^{-2\ell}, x_2^{-2\ell}, \dots, x_{m_{2\ell}}^{-2\ell}, y_1^{-2\ell}, y_2^{-2\ell}, \dots, y_{m_{2\ell}}^{-2\ell}  & & \text{ in degree } -2\ell, \nonumber\\
 	& y_1^{k+i}, y_2^{k+i}, \dots, y_{m_i}^{k+i}  & & \text{ in degree } k+i \text{ for } i=0,1,\dots, 2\ell-1.
 	\end{align}
 	Then we define an element $
 	\omega^0= \sum_{i=0}^{2\ell} \sum_{j=1}^{m_i} d_{dR}x_j^{-i} d_{dR}y_j^{k+i}$ in $ (\Lambda^2\Omega^1_{A})^k$, and  set $ \omega$ to be $ (\omega^0,0,0, \dots) $ as before. Choose an element $H\in A^{k+1}$, the Hamiltonian, satisfying the analogue of \textit{classical master equation}, and define $d$ on $A$ as in Equation (\ref{defn_internal d}) using $H$. We also define the element $\phi \in (\Omega^1_{A})^k$ by the analogue of Equation (\ref{defn_phi}). 
 	 \item[(b)] \cite[Example 5.10]{Brav} When $k/2$ is \emph{odd}, say $k=-4\ell-2$ for $\ell\in \N$, $A$ is freely generated over $A(0)$ by the variables
 	 \begin{align} \label{new local variables for k=-4l-2}
 	 & x_1^{-i}, x_2^{-i}, \dots, x_{m_i}^{-i}   & & \text{ in degree } -i \text{ for } i= 1, 2, \dots, 2\ell, \nonumber \\
 	 & z_1^{-2\ell-1}, z_2^{-2\ell-1}, \dots, z_{m_{2\ell+1}}^{-2\ell-1}  & & \text{ in degree } -2\ell-1, \nonumber\\
 	 & y_1^{k+i}, y_2^{k+i}, \dots, y_{m_i}^{k+i}  & & \text{ in degree } k+i \text{ for } i=0,1,\dots, 2\ell.
 	 \end{align}
 	 We then define an element $
 	 \omega^0= \sum_{i=0}^{2\ell} \sum_{j=1}^{m_i} d_{dR}x_j^{-i} d_{dR}y_j^{k+i}+  \sum_{j=1}^{m_{2\ell+1}} d_{dR}z_j^{-2\ell-1} d_{dR}z_j^{-2\ell-1}$ in $ (\Lambda^2\Omega^1_{A})^k$, and set $ \omega:= (\omega^0,0,0, \dots) $ as before. Choose an element $H\in A^{k+1}$, the Hamiltonian, satisfying the analogue of \textit{classical master equation}
 	 \begin{equation} 
 	 \displaystyle \sum_{i=1}^{2\ell} \sum_{j=1}^{m_i} \dfrac{\partial H}{\partial x_j^{-i}} \dfrac{\partial H}{\partial y_j^{k+i}} + \frac{1}{4} \sum_{j=1}^{m_{2\ell+1}} \Big(\dfrac{\partial H}{\partial z_j^{-2\ell-1}}\Big)^2=0 \text{ in } A^{k+2}.
 	 \end{equation} Define $d$ on $A$ as in Equation( \ref{defn_internal d}) with extra data $ dz_j^{-2\ell-1}:=\dfrac{1}{2}\dfrac{\partial H}{\partial z_j^{-2\ell-1}}. $ 
 	 
 	 Finally, we define the element $\phi \in \Omega^1_{A})^k$ by  \begin{align} \label{defn_phiv2}
 	 \phi&=  \sum_{i=0}^{2\ell} \sum_{j=1}^{m_i} \big[-ix_j^{-i}d_{dR}y_j^{k+i} +(-1)^{i+1}(k+i)y_j^{k+i}d_{dR}x_j^i\big] \nonumber \\
 	 &+k  \sum_{j=1}^{m_{2\ell+1}} z_j^{-2\ell-1} d_{dR}z_j^{-2\ell-1}.
 	 \end{align}
 \end{enumerate}

\begin{observation}\label{observation_vdim symplectic}
  	For $k \not\equiv 2 \mod 4 $, the virtual dimension $\vdim A$ is always \emph{even}.  Otherwise, it can take any value in $ \mathbb{Z} $. More precisely,  for any $k<0$ we have 
 \[\vdim A = 
 \begin{cases*} 
 	 0, & \text{ if }$ k $ \mbox{is odd},\\
 	 \text{even in } \mathbb{Z}, & \text{ if }$ k/2 $ \mbox{is even}, \\
 	 \text{any value in } \mathbb{Z}, & \text{ if }$ k/2 $ \mbox{is odd}. \\
 \end{cases*}
 \]
  \end{observation}
\begin{remark}
The classical Darboux theorem  states that for a symplectic manifold $(X, \omega)$, one can find a local coordinate chart $(U; x_1, \dots, x_n, y_1, \dots, y_n)$ such that $ \omega |_U = \sum_{j} d_{dR}x_j d_{dR}y_j$. Moreover, we can write $\lambda=\sum_{j}x_j d_{dR}y_j$, called the \emph{Liouville form}, such that  $\omega |_U= d_{dR}\lambda $. 

In this derived framework, the element $\phi$ above %can be recognized that it 
	may seem to play the role of $ \lambda. $ However, it is important to notice that $\phi$ is \underline{not} a 1-form (of degree $k$) in the sense of Definition \ref{defn_p form of deg k}, because $d\phi \neq 0.$ %This is in fact the first obstruction to provide an immediate local model for a Darboux-type theorem in the case of shifted contact geometry, and Hence,  one needs to carry out certain algebraic maneuvers to overcome this issue.  In the next section, we shall deal with this problem.
	Therefore, one  needs to modify  $\phi$ to obtain a genuine  1-form of degree $k$.
\end{remark}
%simplification, useful results, proposition 5.7a/proof?

    \section{Shifted contact structures and a Darboux-type theorem} \label{section_Shifted contact structures and Darboux}
\subsection{Basics of classical contact geometry} \label{section_classical contact geo}
  It is very well-known that contact manifolds are viewed as the odd-dimensional analogues of symplectic manifolds. In that respect, they have a number of common features: there is a Darboux theorem providing a local model for such structures; there is no local invariants; and it is more interesting to study their global properties. For more details, we refer to \cite{Geiges}. 
  
  In this section, we shall revisit the basic aspects of contact geometry. There are in fact equivalent ways of describing the notion of a \textit{contact structure}. We prefer to use the  one below. 

\begin{definition} \label{Defn_classical contact strs}
Let $ X $ be a manifold of dimension $ 2n+1 $. \emph{A contact structure} is a smooth field of tangent hyperplanes $ \xi \subset TM$ (of $ \rank 2n $) with the property that for any smooth locally defining 1–form $\alpha$, i.e. $\xi= \ker (\alpha)$, the 2-form $d_{dR}\alpha |_{\xi}$ is non-degenerate. 
%The pair $ (X, \xi) $ is then called a \emph{contact manifold} and $\alpha$ is  called a \emph{local contact form.}% For a point $p\in X$, the pair $(p, \xi_p)$ is called a \emph{contact element.}
\end{definition}
\begin{remark} \label{remark_coorientible and global alpha}
It is also possible to write $\xi= \ker (\alpha)$ with $\alpha$ a \emph{globally} defined contact form on $M$ if and only if $ \xi $ is \emph{coorientable}, by which we mean that the quotient line bundle $TM/\xi$ is trivial. Except some pathological cases, it suffices to work with coorientable contact structures. More details can be found in \cite[\S 1 \& 2 ]{Geiges}.
\end{remark}
Note that %in Definition \ref{Defn_classical contact strs}, by a \emph{non-degenerate 2-form $d_{dR}\alpha |_{\xi}$}, we mean that for each point $p\in X$, the linear map ${d_{dR}\alpha|_{\xi} \ \cdot}: \xi_p \rightarrow \xi_p^*$ defined by $v\mapsto \iota_{v}d_{dR}\alpha$ is an isomorphism. Therefore, 
if $d_{dR}\alpha |_{\xi}$ is non-degenerate, then for each $p\in X,$ $\xi_p = \ker (\alpha_p)$ is a symplectic vector space with a symplectic form $ \omega_p:=d_{dR}\alpha_p |_{\xi_p}.$ Therefore, we also call such 2-form $d_{dR}\alpha |_{\xi}$ \emph{symplectic}.

It follows from the theory of symplectic vector spaces that $\dim \xi_p$ is even and the symplectic form on $\xi_p$ has a canonical form. It means that there exists a symplectic basis $\{e_1, f_1, \dots, e_n, f_n\}$ for $\xi_p$ (and the corresponding dual basis $\{e_1^*, f_1^*, \dots, e_n^*, f_n^*\}$ for $\xi_p^*$) satisfying 
\begin{equation}
\omega_p (e_i, e_j) = 0 =\omega_p (f_i, f_j) \text{ and } \omega_p (e_i, f_j)= -\omega_p (f_j, e_i)=\delta_{ij} \ \ \forall i,j,
\end{equation}so that  $\omega_p$ has the form $\omega_p = \sum_i e_i^* \wedge f_i^*.$ 

Let $ (X, \xi= \ker (\alpha)) $ be a contact manifold of dim $ 2n+1, $ and $p\in X$. Then we have a splitting
\begin{equation} \label{splitting}
T_p X = \xi_p \oplus \ker d_{dR}\alpha_p|_{\xi_p} , 
\end{equation} where $\dim \xi_p=2n$ and $ \dim \ker d_{dR}\alpha_p|_{\xi_p} =1$. In fact, as $ d_{dR}\alpha|_\xi $ is non-degenerate, one can find a local trivialization 
$ \{e_1, f_1, \dots, e_n, f_n, r\} $ of $TM = \ker \alpha \oplus rest$ such that 
\begin{equation} \label{splitting basis}
\ker \alpha = \Span \{e_1, f_1, \dots, e_n, f_n\}  \text{ and } rest=\Span \{r\}.
\end{equation}

Moreover, using this splitting, one can find a unique vector field $ \bf R $, called the \emph{Reeb vector field of $\alpha$}, satisfying $\iota_{\bf R} d_{dR}\alpha = 0$ and $\iota_{\bf R} \alpha = 1.$
\begin{example}
On $\mathbb{R}^{2n+1}$ with cartesian coordinates $(x_1,\cdots,x_n,y_1,\cdots,y_n,z)$, the so-called \emph{standard contact form} is given by
\begin{equation} \label{standard contact form }
\alpha_{std}= \displaystyle -d_{dR}z + \sum_{i=1}^{n}y_i d_{dR}x_i.
\end{equation}Let $\xi \subset T\mathbb{R}^{2n+1}$ be the hyperplane field of $ \rank 2n $ defined by $\alpha_{std},$ i.e. $\xi=\ker \alpha_{std}.$ Then we observe that $$\ker(\alpha_{std}) = \Span \Big \{\dfrac{\partial}{\partial y_j}, y_j\dfrac{\partial}{\partial z}+\dfrac{\partial}{\partial x_j} : j=1,\cdots, n \Big \} \ \text{ and } \  d_{dR}\alpha_{std}= \sum_i d_{dR} x_i \wedge d_{dR}y_i.$$
 Write $A_j={\partial/\partial y_j}$ and $ B_j= y_j{\partial/\partial z}+{\partial/\partial x_j}, $ then it is enough to observe that $d_{dR}\alpha_{std}(A_i, B_j)=\delta_{ij}$ and $d_{dR}\alpha_{std}(A_i,A_j)=0 = d_{dR}\alpha_{std}(B_i,B_j). $ It follows that $ d_{dR}\alpha_{std}|_{\xi} $ is non-degenerate, and hence $ \alpha_{std} $ is a contact form. Moreover, the Reeb vector field of $\alpha_{std}$ is $\textbf{R} = {\partial / \partial z}.$

\end{example}

As in the symplectic case, there is a Darboux-type theorem for contact structures. It basically says that all contact structures can be locally given as in Equation (\ref{standard contact form }). More formally, we have

\begin{theorem}[Darboux Theorem for contact structures]
Let $(X,\alpha)$ be a contact manifold of dimension $2n+1$, and $p \in X$. Then there exists a local coordinate system $(U; x_1,\cdots,x_n,y_1,\cdots,y_n,z)$ around $p$ such that $p=(0,0,\cdots,0)$ and
\begin{equation} 
\alpha|_U= \displaystyle -d_{dR}z + \sum_{i=1}^{n}y_i d_{dR}x_i.
\end{equation}
\end{theorem}

\subsection{Shifted contact structures and Darboux models} \label{section_shifted contact structures and Darboux forms}
	In this section, we  provide an appropriate analogue of Definition \ref{Defn_classical contact strs} for  derived spaces %$\K$-schemes (or more generally, for derived stacks) 
	and study the local theory of shifted contact structures.  	
	Let $\bf X$ be a locally finitely presented derived (Artin) stack. Then we have:
\begin{definition}\label{defn_preshiftedcontact}
 A \bfem{pre-$ k $-shifted contact structure} on $\bf X$ consists of a perfect complex $\mathcal{K}$ on $\bf X$ with a morphism $\kappa: \mathcal{K} \rightarrow \mathbb{T}_{{\bf X}}$ of perfect complexes whose cone $Cone(\kappa)$ is of the form $L[k]$, up to quasi-isomorphism, where $L$ is a line bundle\footnote{In the spirit of Remark \ref{remark_coorientible and global alpha}, we say that a pre-$ k $-shifted contact structure is \emph{coorientable} if $L$ in the data is trivial.}. Denote such a structure on $\bf X$  by $(\mathcal{K}, \kappa, L)$.

\end{definition}

\begin{definition}\label{defn_shiftedcontact}
	We say that a pre-$ k $-shifted contact structure $(\mathcal{K}, \kappa, L)$ on $\bf X$ is a \bfem{$ k $-shifted contact structure } if locally on $\bf X$, where $L$ is trivial, the induced $k$-shifted 1-form $\alpha: \mathbb{T}_X \rightarrow \mathcal{O}_X[k]$ is such that the induced map $d_{dR}\alpha|_{\mathcal{K}} \ \cdot:= \kappa^{\vee}[k] \circ( d_{dR}\alpha \ \cdot) \circ \kappa : \mathcal{K}\rightarrow\mathcal{K}^{\vee}[k]$ is a weak equivalence. In that situation, we say the $ k $-shifted 2-form $d_{dR}\alpha$ is \emph{non-degenerate on $\mathcal{K}$}.  
	
	Moreover, when $k\leq 0$, the fiber-cofiber sequence $\mathcal{K} \rightarrow \mathbb{T}_{{\bf X}} \rightarrow L[k]$ in $QCoh(X)$ splits  Zariski
	locally for any derived scheme. Call such local form a \bfem{$k$-contact form}.%\footnote{For the coorientable shifted contact structures, the induced form and the splitting will be globally defined. }.
\end{definition}

Let $\bf X$ be a locally finitely presented  derived Artin stack with a $k$-shifted contact structure $(\mathcal{K}, \kappa, L)$. Recall from Yoneda's lemma, ${\bf X}(A) \simeq Map_{dPstk}(\spec A, {\bf X})$, and hence any $A$-point $p\in {\bf X}(A) $ can be seen as a morphism $p: \spec A \rightarrow \bf X$ of derived pre-stacks. Then, let us consider the pair $ (p, \alpha_p )$, with $  p \in {\bf X} (A), \ \alpha_p \in p^*(\mathbb{L}_{{\bf X}}[k])$ a  $k$-contact form. For $A \in cdga_{\K}$, there is a $\mathbb{G}_m(A)$-action on the pair $ (p, \alpha_p )$ by 
\begin{equation*}
f \triangleleft (p, \alpha_p):= (p, f \cdot \alpha_p).
\end{equation*} Denote by $H^0$  the functor sending $A \mapsto H^0(A).$ Denote the image under $H^0$ of an element $f$ simply by $f^0.$ Note that localizing $ A $ if necessary, w.l.o.g. we may assume that the image $f^0$ is always invertible. It follows that $f^0$ lies in $(A^0)^{\times}$, which is by definition $\mathbb{G}_m(A^0)=(A^0)^{\times}$.

\begin{observation} \label{observation_cocones}
	If ${\bf X}, (p, \alpha_p )$, and the $\mathbb{G}_m(A)$-action are as above, then for an element $f\in \mathbb{G}_m(A)$, we can obtain $Cocone(f\cdot \alpha_p) \simeq Cocone(\alpha_p)$ by using the invertibility of $f.$ 
\end{observation}

Because our results in this paper are all about the \emph{local structure theory} of contact derived stacks, we will  focus on the \emph{refined affine case}  in the sense of Theorem \ref{localmodelthm}. In this regard, we will always assume that our \textit{refined local models} are given in terms of minimal standard form cdgas.

 From Proposition \ref{proposition_L as a complex of H^0 modules}, on a refined affine neighborhood, say on $\spec A$ with $A$ a minimal standard form cdga, the perfect complexes $\mathbb{T}_A, \mathbb{L}_A$, when restricted to ${\spec H^0(A)}$, are both free finite complexes of $H^0(A)$-modules. In that case, Definitions \ref{defn_preshiftedcontact} and \ref{defn_shiftedcontact}, and Observation \ref{observation_cocones} will reduce to the following local descriptions, where  $\mathcal{K}, Cone(\kappa)$ will be equivalent to the ordinary $\ker (\alpha), \coker(\kappa) $, respectively in $D(Mod_A)$;  and  $L$ in the splitting corresponds to the  line bundle generated by the Reeb vector field of the classical case.   
 
 \paragraph{Shifted contact structures with (nice) local models.}%\label{defn_shiftedcontact locally}
 Recall that the \textit{(mapping) cone} of a morphism $f:A\rightarrow B$ in some homotopical category is a realization of the homotopy cofiber of $f$. That is, it is the homotopy pushout satisfying universally the homotopy diagram 
 
 \begin{equation} \label{eqn_homotopy cokernel}
 	\begin{tikzpicture}
 		\matrix (m) [matrix of math nodes,row sep=2em,column sep=3 em,minimum width=1.5 em] {
 			A   & \star  \\
 			B &  Cone(f). \\};
 		\path[-stealth]
 		(m-1-1) edge  node [left] {{\small $ f $} } (m-2-1)
 		(m-1-1) edge  node [above] { } (m-1-2)
 		%(m-1-1) edge  node [below] {} node [below] {{\small stacks}} (m-2-2)
 		%(m-1-1) edge  node [below] {} node [below] {{\small  higher stacks}} (m-3-2)
 		(m-2-1) edge  node  [below] { } (m-2-2)
 		
 		(m-1-2) edge  node [right] { } (m-2-2);
 		%edge [dashed,-] (m-2-1);
 	\end{tikzpicture}
 \end{equation} It is also called the \textit{homotopy cokernel} of $f$ or the \textit{weak quotient} of $B$ by the image of $A$ under $f$.
 
 As a dual notion, the \emph{(mapping) cocone} of a morphism $f:A\rightarrow B$, with $B$ a pointed object, in a homotopical category is a particular realization of the homotopy fiber of $f$ (i.e. of the homotopy pullback of the point along $f$). In that case, the following diagram homotopy commutes:
 \begin{equation} \label{eqn_homotopy kernel}
 	\begin{tikzpicture}
 		\matrix (m) [matrix of math nodes,row sep=1.5em,column sep=3 em,minimum width=1.5 em] {
 			Cocone(f)   & \star  \\
 			A &  B. \\};
 		\path[-stealth]
 		(m-1-1) edge  node [left] { } (m-2-1)
 		(m-1-1) edge  node [above] { } (m-1-2)
 		%(m-1-1) edge  node [below] {} node [below] {{\small stacks}} (m-2-2)
 		%(m-1-1) edge  node [below] {} node [below] {{\small  higher stacks}} (m-3-2)
 		(m-2-1) edge  node  [below] {{\small $ f $} } (m-2-2)
 		
 		(m-1-2) edge  node [right] { } (m-2-2);
 		%edge [dashed,-] (m-2-1);
 	\end{tikzpicture}
\end{equation} It is also called the \textit{homotopy kernel} of $f$.

\begin{remark} \label{remark_cocone/cone vs kernel/cokernel}
	In the case of a morphism (perfect) of complexes $f:A\rightarrow B$, we have $\star=0$ (as the initial and terminal objects) and that the complex $\ker(f)$, which is the subcomplex of $A$ formed by the (strict) kernels $\{ \ker (f_n) \},$  commutes the diagram
	\begin{equation} 
		\begin{tikzpicture}
			\matrix (m) [matrix of math nodes,row sep=2em,column sep=3 em,minimum width=1.5 em] {
			\ker(f)   & \star  \\
				A &  B, \\};
			\path[-stealth]
			(m-1-1) edge  node [left] { $ i $ } (m-2-1)
			(m-1-1) edge  node [above] { } (m-1-2)
			%(m-1-1) edge  node [below] {} node [below] {{\small stacks}} (m-2-2)
			%(m-1-1) edge  node [below] {} node [below] {{\small  higher stacks}} (m-3-2)
			(m-2-1) edge  node  [below] {{\small $ f $} } (m-2-2)
			
			(m-1-2) edge  node [right] { 0 } (m-2-2);
			%edge [dashed,-] (m-2-1);
		\end{tikzpicture}
	\end{equation} which essentially means $ f\circ i = 0 $. This is in fact a strict pullback diagram. From the universality of the homotopy kernel of $f$, there is a natural inclusion $\ker (f) \hookrightarrow Cocone(f)$ from the strict fiber to the homotopy fiber (i.e. from the strict kernel to the homotopy kernel).  %the complex $\ker (f)$ is in fact equivalent to $cocone(f)$.
This map is in fact an equivalence if $f$ is a fibration (i.e. a surjective morphism of complexes).

$\coker (f)$, on the other hand, is the quotient complex of $B$	formed by the cokernels $\{ \coker (f_n) \}.$ By definition, $\coker (f)$ commutes the diagram
	\begin{equation} 
		\begin{tikzpicture}
			\matrix (m) [matrix of math nodes,row sep=2em,column sep=3 em,minimum width=1.5 em] {
				A   & \star  \\
				B &  \coker (f). \\};
			\path[-stealth]
			(m-1-1) edge  node [left] {{\small $ f $} } (m-2-1)
			(m-1-1) edge  node [above] { } (m-1-2)
			%(m-1-1) edge  node [below] {} node [below] {{\small stacks}} (m-2-2)
			%(m-1-1) edge  node [below] {} node [below] {{\small  higher stacks}} (m-3-2)
			(m-2-1) edge  node  [below] { } (m-2-2)
			
			(m-1-2) edge  node [right] { 0 } (m-2-2);
			%edge [dashed,-] (m-2-1);
		\end{tikzpicture}
	\end{equation} From the universality of the homotopy cokernel of $f$, there is a natural map $Cone(f)\rightarrow \coker (f)$ from the homotopy cofiber to the strict cofiber (i.e. from the homotopy cokernel to the strict cokernel). This map is in fact an equivalence if $f$ is a cofibration (i.e. an injective morphism of complexes).

\end{remark}
%\newpage
Now, in the presence of Remark \ref{remark_cocone/cone vs kernel/cokernel}, we revisit Definitions \ref{defn_preshiftedcontact} and \ref{defn_shiftedcontact} with (nice) local models.	Let $\textbf{X}=\spec A$ be an affine derived $\mathbb{K}$-scheme for $A$ a minimal standard form cdga, endowed with a $k$-shifted contact structure $(\mathcal{K}, \kappa, L)$ for $k<0$. Assume also that $L$ is trivial on $\spec A.$
From definitions, the triangle $\mathcal{K} \rightarrow \mathbb{T}_{{\bf X}} \rightarrow L[k]$ splits and $ Cone (\kappa) $ is of the form a line bundle $L[k]$. In that case, we have the homotopy commuting square 
\begin{equation} 
	\begin{tikzpicture}
		\matrix (m) [matrix of math nodes,row sep=2em,column sep=3 em,minimum width=1.5 em] {
			\mathcal{K}   & \star  \\
			\mathbb{T}_{{\bf X}} &  Cone (\kappa). \\};
		\path[-stealth]
		(m-1-1) edge  node [left] {{\small $ \kappa $} } (m-2-1)
		(m-1-1) edge  node [above] { } (m-1-2)
		%(m-1-1) edge  node [below] {} node [below] {{\small stacks}} (m-2-2)
		%(m-1-1) edge  node [below] {} node [below] {{\small  higher stacks}} (m-3-2)
		(m-2-1) edge  node  [below] { } (m-2-2)
		
		(m-1-2) edge  node [right] { 0 } (m-2-2);
		%edge [dashed,-] (m-2-1);
	\end{tikzpicture}
\end{equation} If  $ \kappa: \mathcal{K} \rightarrow \mathbb{T}_{{\bf X}} $ is a monomorphism, the homotopy cofibers are equivalent to strict cofibers by Remark \ref{remark_cocone/cone vs kernel/cokernel}. Thus, $ Cone (\kappa) $ is equivalent to the ordinary $\coker(\kappa)$.

On $\textbf{X}=\spec A$, where $L$ is trivial, $L[k]$ is  of the form $\mathcal{O}[k].$ Recall that since $A$ is a minimal standard form cdga, we have $\mathbb{L}_{{\bf X}}\simeq \mathbb{L}_A$ in $D(Mod_A)$ and the perfect complexes $\mathbb{T}_A, \mathbb{L}_A$, when restricted to ${\spec H^0(A)}$, are both free finite complexes of $H^0(A)$-modules.  Let $\alpha$ be a  $ k$-contact form, with the underlying $k$-shifted 1-form $\alpha: \mathbb{T}_A \rightarrow A [k]$, then we have
\begin{equation} 
	\begin{tikzpicture}
		\matrix (m) [matrix of math nodes,row sep=2em,column sep=3 em,minimum width=1.5 em] {
			Cocone(\underline{\alpha}) &                &       \\
			               & \mathcal{K}   & \star  \\
							& \mathbb{T}_A &  Cone (\kappa) \simeq \mathcal{O}[k]. \\
						};
		\path[-stealth]
		%\path[->]
		(m-2-2) edge  node [left] {{\small $ \kappa $} } (m-3-2)
		(m-2-2) edge  node [above] { } (m-2-3)
		(m-1-1) edge  node [below] {{\small $ \simeq $}} (m-2-2)
		(m-1-1) edge [bend right=25,looseness=1]  (m-3-2)
		(m-1-1) edge [bend left=15,looseness=1]  (m-2-3)
		(m-3-2) edge  node  [below] { } (m-3-3)
		
		(m-2-3) edge  node [right] { 0 } (m-3-3)
		(m-3-2) edge  [bend right=30,looseness=1.1] node [below] {$\alpha$ } (m-3-3);
		
		%edge [dashed,-] (m-2-1);
	\end{tikzpicture}
\end{equation} Here both outer and inner squares (resp. the homotopy fiber and the homotopy cofiber)  commute, and $\mathcal{K}$ is  equivalent to $Cocone(\underline{\alpha})$, with $\underline{\alpha} \in \mathcal{A}^1(\spec A, k)$ a map\footnote{The map from the epi-mono factorization of $\alpha: \mathbb{T}_A \twoheadrightarrow \Ima \alpha \hookrightarrow \mathcal{O}[k]$, where $ \underline{\alpha}:\mathbb{T}_A \rightarrow \Ima \alpha. $} of perfect complexes. From Remark \ref{remark_cocone/cone vs kernel/cokernel}, there is a natural map $\ker (\underline{\alpha}) \hookrightarrow \mathcal{K}$ from the strict kernel to the homotopy kernel of $\underline{\alpha}$. Since $ \underline{\alpha} $ is an epimorphism, the natural map $\ker (\underline{\alpha}) \hookrightarrow \mathcal{K}$ is then an equivalence. Thus, on suitable local models, we can use the strict kernel to represent the homotopy kernel.

 Now, we consider an explicit general form of the underlying $k$-shifted 1-form $\alpha: \mathbb{T}_A \rightarrow A [k]$ on $\spec A$ with $A$ as above.  Note that it is the minimal (at $p\in \spec H^0(A)$) compared to all other cdgas quasi-isomorphic to $A$ (at $p\in \spec H^0(A)$). Explicitly, Letting $A=A(-k)$, $A$ is  the free graded  algebra over $A(0)$ generated by the variables $ x_1^{-i}, x_2^{-i}, \dots, x_{m_i}^{-i}$, with $m_i \in \mathbb{Z}$, for $i=1, \dots, -k$ such that $d_{dR}x_1^{-i}, d_{dR}x_2^{-i}, \dots, d_{dR}x_{m_i}^{-i}, \ i=1, \dots, -k,$ is a $A$-basis for $ \Omega_A^1 $. Write $\underline{\alpha}=\alpha= \sum_{i,j} \alpha_j^{k+i}d_{dR}x_j^{-i}$ with $\alpha_j^n\in A^n$. Notice that by definition, we have $\ker (\alpha)=\ker(\underline{\alpha})$. Thus, when we consider the strict kernel on such local models, we simply write  $\ker (\alpha).$

Therefore, the discussion above implies that the data of \emph{(1) a submodule $\mathcal{K}$ of $\mathbb{T}_A=Der(A)$ with the natural inclusion $i:\mathcal{K}\hookrightarrow Der(A)=(\Omega^1_{A})^{\vee}$}; \emph{(2) a morhism $\mathbb{T}_A \rightarrow L[k]\simeq\coker (i)$ for $ L[k] $ a shifted line bundle;} and \emph{(3) the non-degeneracy condition on the induced map $\mathbb{T}_A \rightarrow \mathcal{O}[k]$} will define a shifted contact structure on $\textbf{X}=\spec A$, where $A$ is a  (minimal) standard form cdga. That is, we  obtain the following proposition/definition for nice affine derived schemes:

\bfem{Definition with nice local models.} For a (minimal) standard form cdga $A$ and $k<0$, a \bfem{$ k $-shifted (strict) contact structure on $\spec A$} is a submodule $\mathcal{K}$ of $\mathbb{T}_A=Der(A)$ with the natural inclusion $i:\mathcal{K}\hookrightarrow Der(A)=(\Omega^1_{A})^{\vee}$ such that $\mathcal{K}\simeq \ker (\alpha)$ for a $k$-shifted 1-form $\alpha$    with the property that the $k$-shifted 2-form $d_{dR}\alpha$ is non-degenerate  on $\ker (\alpha)$ and the complex $Cone(i: \mathcal{K} \hookrightarrow \mathbb{T}_A)\simeq \coker(i)$ is the quotient complex and of the form $L[k]$, with $L$ a line bundle.
%\end{definition}

\begin{observation} \label{observation_equivalent contact strs}
	Let $\textbf{X}=\spec A$ be an affine derived $\mathbb{K}$-scheme for $A$ a minimal standard form cdga. For any $f\neq 0$ in $A$ and any $k$-shifted contact form $\alpha$, one has $\ker (\alpha) \simeq \ker (f\alpha).$ Hence, both define equivalent contact structures on $\bf X$. In fact, this follows from the fact that the contraction operation $\iota_{Y}$ on the de Rham algebra $DR(A)$ with a homogeneous vector field $Y$ is the unique derivation of  degree $|Y|+1$ such that $\iota_{Y}g=0$ and $\iota_{Y}d_{dR}g=Y(g)$ for all $g\in A.$ Therefore,
	\begin{equation*}
	\iota_{Y} (f \alpha) = (\iota_{Y}f)\cdot \alpha + (-1)^{|Y|+1}f\cdot \iota_Y \alpha = (-1)^{|Y|+1}f\cdot \iota_Y \alpha.
	\end{equation*} %As $f\neq 0$, $ \iota_Y \alpha =0 \iff \iota_{Y} (f \alpha)=0. $ 
	Adopting the classical terminology (in terms of non-integrable distributions), on refined local models, we sometimes call the subcomplex $\ker (\alpha)$ of $\mathbb{T}_{{\bf X}}$ a \textit{contact structure} and the corresponding $ k$-shifted 1-form $ \alpha $ a \textit{(locally) defining $k$-contact form}.
\end{observation}
In what follows, we give a prototype construction for  $k$-shifted contact forms, which is similar to the previous case of shifted symplectic Darboux models. In brief, we will coherently extend the symplectic case studied in \cite[Example 5.8]{Brav} using similar constructions.
\begin{example} \label{model example}
	In this example, fixing $\ell\in \N$, we construct a standard form cdga $ A = A(n)$ with the induced $k$-shifted contact structure on $\spec A$, for $k=-2\ell -1=-n$.
	
	\bfem{Step-1: Construction of a cdga $(A,d)$.} First, we consider a smooth $\K$-algebra $A(0)$ of dimension $m_0$. We assume that there exist degree 0 variables $x_1^0, x_2^0, \dots, x_{m_0}^0$ in $A(0)$ such that $d_{dR}x_1^0,\dots,d_{dR}x_{m_0}^0$ form a basis for $\Omega_{A(0)}^1$ over $A(0).$ This choice can be made by localizing $A(0)$ if necessary.	
	Next, choosing non-negative integers $m_1,\dots, m_{\ell}$, define  a \bfem{commutative graded algebra} $A$ to be the free graded $\K$-algebra over $A(0)$ generated by the variables  
	\begin{align} \label{var set1}
		& x_1^{-i}, x_2^{-i}, \dots, x_{m_i}^{-i}& &\text{ in degree } -i \ \ \ \text{ for } i= 1,  \dots, \ell, \nonumber \\
		& y_1^{k+i}, y_2^{k+i}, \dots, y_{m_i}^{k+i}& & \text{ in degree } k+i \ \text{ for } i=1,\dots, \ell, \nonumber \\
		& z^k, y_1^{k}, y_2^{k}, \dots, y_{m_0}^{k}& & \text{ in degree } k,
	\end{align} where we call $z^k$  the \emph{distinguished variable (of $\deg k$)}. It follows that  $\Omega^1_{A}$ is  the free  $A$-module of finite rank with an $A$-basis $ \big \{d_{dR}x_j^{-i}, d_{dR}y_j^{k+i}, d_{dR}z^k : i= 0, 1, \dots, \ell, \ j= 1,2, \dots, m_i\big \}.$

Choose an element $H\in A^{k+1}$ satisfying  the classical master equation (\ref{defn_CME}). Here, due to  degree reasons, $H$ does not involve  any of $z^k,y^k_j$'s.
Then we define the \bfem{internal differential} $d$ on $A$   by the equations \begin{align} \label{defn_internal d contact}
	d|_{A(0)}&=0; \ dx_j^{-i} =  \dfrac{\partial H}{\partial y_j^{k+i}} \text{ for all } i>0,j; \  \ dy_j^{k+i} =  \dfrac{\partial H}{\partial x_j^{-i}} \text{ for all } i,j; \text{ and } \nonumber \\ -kdz^k&= H+d\Big[\sum_{i,j} (-1)^{i} ix_j^{-i} y_j^{k+i} \Big].
\end{align}  %Moreover, we simply have $dH=0.$

Here, the condition on $H$ implies that $d^2=0$ on each generator \cite{Brav}. By construction, $A$ is then a standard form cdga with $A=A(n=2\ell+1)$ which is defined inductively by adjoining free modules $M^{-i}=\langle x_1^{-i}, x_2^{-i}, \dots, x_{m_i}^{-i} \rangle_{A(i-1)}$ for $ i= 1, 2, \dots, \ell $; $M^{k+i}=\langle y_1^{k+i}, y_2^{k+i}, \dots, y_{m_i}^{k+i} \rangle_{A(-k-i-1)}$ for $ i= 1, \dots, \ell; $ and $M^{k}=\langle y_1^{k}, y_2^{k}, \dots, y_{m_0}^{k}, z^k \rangle_{A(-k-1)}$. Also, the cdga $A$ has  $\vdim(A)=-1$.

	\bfem{Step-2: Pre-contact data.} %localizing $A$ at $p\in \cspec H^0(A)$ if necessary and using minimality of $A$ at $p$, 
 Define an element $\alpha \in (\Omega_A^1)^k$ by %called the primitive \emph{$k$-contact element}, by
\begin{equation} \label{the local primitive contact model}
	\alpha= d_{dR}z^k+ \displaystyle \sum_{i=0}^{\ell} \sum_{j=1}^{m_i} y_j^{k+i}d_{dR}x_j^{-i}.
\end{equation} Let us first show $d\alpha=0.$ Start with the computation $ d_{dR}\alpha= \sum_{i,j} d_{dR}x_j^{-i} d_{dR}y_j^{k+i} \in  (\Lambda^2\Omega^1_{A})[k] $. \\
\cite[Example 5.8]{Brav} proves that      $d_{dR}\alpha$ is both $d$- and $d_{dR}$-closed   such that for the choices\footnote{Alternatively, one may take $\phi= k\sum_{i,j} x_j^{-i}d_{dR}y_j^{k+i}$ or $ \phi=  k\sum_{i=0}^{\ell} \sum_{j=1}^{m_i} y_j^{k+i}d_{dR}x_j^{-i}$ replacing $H, \phi$ by suitable $H+d[\cdots]$ and $\phi+\dR [\cdots]$.} of elements $$\phi:= \sum_{i=0}^{\ell} \sum_{j=1}^{m_i} [(-i)x_j^{-i}d_{dR}y_j^{k+i}+(k+i)y_j^{k+i}d_{dR}x_j^{-i}]\in (\Omega^1_A)^k$$ and $H \in A^{k+1}$ as above, the pair $(\phi,H)$ is a solution to the equations \begin{equation} \label{important relations}
dH=0 \text{ in } A^{k+2}, \ d_{dR}H+d\phi=0 \text{ in } (\Omega^1_{A})^{k+1}, \text{ and } d_{dR}\phi = kd_{dR}\alpha.
\end{equation} Observe that  $\phi  + \dR \big[ \sum_{i,j} (-1)^{i} ix_j^{-i} y_j^{k+i} \big]=  k\sum_{i,j} y_j^{k+i}d_{dR}x_j^{-i}$, and hence we can write 
\begin{align} \label{local contact model}
	k\alpha&=kd_{dR}z^k+ \displaystyle k\sum_{i=0}^{\ell} \sum_{j=1}^{m_i} y_j^{k+i}d_{dR}x_j^{-i}=kd_{dR} z^k + \phi + \dR[\cdots].
\end{align} Here $k\alpha \in \Omega^1_{A}[k]$ is a representative of $\phi$, but with  $d\phi=-d_{dR} H$.   
As we let ${-kdz^k=H+d[\cdots]}$ in (\ref{defn_internal d contact}), we obtain $d(k\alpha)=-k\dR \circ dz^k + d\phi + d \circ \dR [\cdots] = \dR(H+d[\cdots])-\dR H - \dR \circ d [\cdots]=0$ using (\ref{important relations}). In that case, $\alpha$ is then \emph{$d$-closed, and hence a 1-form of degree $k.$} In other words, the pair $(H,k\alpha-kd_{dR}z^k-\dR[\cdots])$ satisfies the  equations (\ref{important relations}), which implies $d\alpha=0.$ 

 We now show that there is a natural pre-$k$-shifted contact structure on $\spec A$ defined by $\ker \alpha.$ For $i=0,\dots,\ell, \text{ and } 1\leq j\leq m_i$, we denote the vector fields annihilating $\alpha$ by
	  \begin{align}
	  	%\sigma_{ij}^{0} &=ix^{-i}_j\partial/\partial x_j^{-i} + (k+i)y^{k+i}_j\partial/\partial y^{k+i}_j & & \text{ in degree }  0, \nonumber \\
	  	\zeta_j^{i}&= \partial/\partial x_j^{-i}-y^{k+i}_j\partial/\partial z^{k} & & \text{ in degree } i, \nonumber \\
	  	\eta^{-k-i}_j &= \partial/\partial y^{k+i}_j & & \text{ in degree }  -k-i.
	  \end{align} %Notice that $ \sigma_{ij}^0$ can also be written by using the other vector fields\footnote{Observe that $ k\sigma_{ij}^{0}= ix^{-i}_j	\zeta_j^{i}+ (k+i)y^{k+i}_j\eta^{-k-i}_j$.} $ 	\zeta_j^{i} \text{ and }  \eta^{-k-i}_j$. 
  We thus obtain $\ker \alpha= \Span_{A} \{ \zeta_{j}^{i}, \ \eta^{-k-i}_{j} \} \hookrightarrow  \mathbb{T}_A $ over $\spec A$, with the  $\mathrm{Tor}$-amplitude $[0,-k]$.  When restricted to $ \cspec H^0(A)$,    $\ \faktor{\mathbb{T}_A}{\ker \alpha}$ is then generated by the  vector field $\partial/\partial z^{k}$ as an $H^0(A)$-module, and hence it is equivalent to the complex concentrated in degree $-k$. Then the quotient, denoted by $Rest$,  can also be identified with a $k$-shifted line bundle  on  $\spec A$ such that we get the  splitting 
	  \begin{equation} 
	  \mathbb{T}_A|_{\cspec H^0(A)}=  \ker \alpha |_{\cspec H^0(A)} \oplus Rest|_{\cspec H^0(A)},
	  \end{equation}where $ \mathbb{T}_A|_{\cspec H^0(A)} $ is a complex of free $ H^0(A)$-modules by Proposition \ref{proposition_L as a complex of H^0 modules}, %such that   for each homogeneous degree $n$,   we have $ (\mathbb{T}_A|_{\cspec H^0(A)})^n=  (\ker \alpha |_{\cspec H^0(A)})^n \oplus Rest^n|_{\cspec H^0(A)}$ 
  and
	  \begin{align*}
	  \ker \alpha|_{\cspec H^0(A)}&= \langle   \zeta_j^{i}, \ \eta^{-k-i}_j : \  i=0,1,\dots \ell, \ j=1,\dots, m_i\rangle_{A(0)}, \\
	  Rest|_{\cspec H^0(A)} &=\langle \partial/\partial z^k \rangle_{A(0)} . 
	  \end{align*}
  Hence, we get a \bfem{pre-$k$-shifted contact structure} on $\spec A$ induced by the $k$-shifted 1-form $\alpha$ in the sense of Definition \ref{defn_preshiftedcontact}.  
  Call such $\alpha$  a \bfem{pre-$k$-contact form}.
 
	 \bfem{Step-3: From pre-contact to contact.} So far, we have obtained a pre-$k$-contact form $\alpha \in (\Omega_A^1)^k$ as in (\ref{the local primitive contact model}) using  the  variables (\ref{var set1}), the differential $d$ given by (\ref{defn_internal d contact}), and the  equations in (\ref{important relations}). It remains to show that $d_{dR}\alpha$ is non-degenerate on $ {\ker \alpha}$.
	To this end,  it suffices to prove the non-degeneracy of the map  $d_{dR}\alpha|_{\ker \alpha}\otimes_A \mathrm{id}_{H^0(A)}: \ker \alpha \otimes_A H^0(A) \rightarrow (\ker \alpha)^{\vee}[k] \otimes_A H^0(A). $
	
	We first observe that, at $p\in \cspec H^0(A)$, $d_{dR}\alpha|_{p}$ maps \begin{align*}
		&\langle \partial/\partial x_1^{-i}|_{p}, \dots, \partial/\partial x_{m_i}^{-i}|_{p} \rangle_{\K} \xrightarrow{\sim} \langle d_{dR}y_1^{k+i}|_{p},\dots, d_{dR}y_{m_i}^{k+i}|_{p} \rangle_{\K}, \\ &\langle \partial/\partial y^{k+i}_1|_{p},\dots,\partial/\partial y^{k+i}_{m_i}|_{p} \rangle_{\K} \xrightarrow{\sim} \langle d_{dR}x_1^{-i}|_{p},\dots, d_{dR}x_{m_i}^{-i}|_{p} \rangle_{\K}
	\end{align*} isomorphically for all $i$, with $d_{dR} \alpha (-, \partial/\partial z^k|_{p})=0$. It follows that,  at $p\in \spec H^0(A)$, we have the identifications 
	\begin{equation}\label{identification T_A/ker d_DR alpha}
		\Big(\faktor{\mathbb{T}_A}{\ker d_{dR}\alpha}\Big)|_p  \simeq \langle \partial/\partial x_j^{-i}|_p, \partial/\partial y^{k+i}_j|_p : \forall i,j\rangle_{\K}		\simeq (\ker \alpha)|_p.
	\end{equation} 
	Thus, the maps $(d_{dR}\alpha|_{\ker \alpha})^{i}|_p: (\ker \alpha|_p)^{i}\rightarrow (\ker^{\vee} \alpha|_p)^{k+i} $ are isomorphisms at $p$, hence in a neighborhood of $p$. So, localizing $A$ at $p$ if necessary, the map $d_{dR}\alpha|_{\ker \alpha}\otimes_A \mathrm{id}_{H^0(A)}$ is an isomorphism of complexes, and hence a quasi-isomorphism. This proves the non-degeneracy. 
	Thus, we get a \bfem{$k$-shifted contact structure} on $\spec A$  with the \bfem{$k$-contact form $\alpha$} in the sense of Definition \ref{defn_shiftedcontact}.
	\begin{definition}
		If  $A, d, \alpha$  are as above, %would serve as the desired contact model. 
		we then say $A, \alpha$  are \bfem{in contact Darboux form}. 
	\end{definition}
\end{example}

\begin{remark} \label{remark_alternative and general versions of alpha}
\begin{enumerate}
\item Let $k, A, H,\phi$ be as in Example \ref{model example}. Alternatively, we can define the \bfem{differential} $d$ on $A$ as in (\ref{defn_internal d contact}), but with $-kdz^k=H$, instead. In that case, we then define the element $ \alpha' \in \Omega_A^1[k]$  by \begin{align*}
	\alpha'&= d_{dR}z^k+  \sum_{i,j} \Big[-\dfrac{i}{k} \ x_j^{-i}d_{dR}y_j^{k+i}+\frac{k+i}{k} \ y_j^{k+i}d_{dR}x_j^{-i}\Big] =d_{dR}z^k+\phi/k.
\end{align*}Observe that  $d\alpha'=  0$ as well due to the new choice of $dz^k$. Modifying Example \ref{model example} accordingly shows that such element $\alpha'$ also induces a $k$-contact structure on $\spec A$, with the contact data determined similarly. Details are left to the reader. 

\item Note  that the general expressions like $``d_{dR} z^k + \phi/k +\dR [\cdots]/k"$ in Example \ref{model example} will still be valid for the other cases $ (a) \ k\equiv0 \mod 4,$  and  $(b) \ k\equiv2 \mod 4  $. Equations (\ref{new local variables for k=-4l}) $-$ (\ref{defn_phiv2}) show that the other cases in fact involve modified versions of  $H, d, $ and  $\phi$ with some possible extra terms. In any case, the modified $A, \alpha$ would also serve as the desired contact model.
\end{enumerate}
\end{remark}  Following the same terminology, we refer to  all possible models mentioned in Remark \ref{remark_alternative and general versions of alpha} and Example \ref{model example} as \bfem{(standard) contact Darboux forms}. These cases also suggest  that suitable modifications using $H+d[\cdots]$ and $\phi+\dR [\cdots]$ may  lead to alternative versions of such forms.

\subsection{A Darboux-type theorem for negatively shifted contact derived schemes} \label{section_darboux theorem}
In what follows, we give the proof of Theorem \ref{THM1}, which essentially says that for $k<0$, every $k$-shifted contact derived $\K$-scheme $\bf X$ is locally equivalent to $(\spec A, \alpha_0)$ for  a minimal standard form cdga $A$ and a $k$-contact form $\alpha_0$. % similar to the ones in Example \ref{model example} (with some minor modifications, if necessary).
 Our result is a ``contact" variation of \cite[Theorem 5.18.]{Brav}. More precisely, we have:

\begin{theorem} \label{contact darboux}
	Let $\bf X$ be a (locally finitely presented) derived $\mathbb{K}$-scheme with a $k$-shifted contact structure $(\mathcal{K}, \kappa, L)$ for $k<0$, and $x\in \bf X$. Then there is a local contact model $\big(A,  \alpha_0 \big)$  and $p \in \cspec H^0(A)$ such that $i: \spec A \hookrightarrow \bf X$ is an open inclusion with $i(p)=x$, \ $A$ is a standard form cdga that is minimal at  $p$, and $\alpha_0$ is a $k$-shifted contact form on $\spec A$ such that $A, \alpha_0$ are in standard contact Darboux form.%, and $i^*(\alpha) \sim \alpha_0$ in the space of $k$-shifted 1-forms. 

\end{theorem} Note that for $k<0$ odd,  such pair $(A, \alpha_0)$ is explicitly described in Example \ref{model example}, see Equations (\ref{var set1})$, $ (\ref{the local primitive contact model}),  (\ref{local contact model}), and Remark \ref{remark_alternative and general versions of alpha}. For the other cases, one should use another sets of variables as in Equations (\ref{new local variables for k=-4l}) and (\ref{new local variables for k=-4l-2}), and modify $H, \phi, d$ accordingly.

Before giving the proof, we begin by some remarks and simplifying assumptions.
%\paragraph{Proof of Theorem \ref{contact darboux}.}

 \begin{remark} \label{remark_setup before theproof1}
We first note that as $\bf X$ is a locally finitely presented derived $\mathbb{K}$-scheme, there exists a cover by affine  derived $\mathbb{K}$-subschemes of finite presentation. Thus, for $x\in \bf X$ we can choose a cdga $B$ of finite presentation with a Zariski open inclusion $ \iota: \spec B \hookrightarrow \bf X $ and a unique $q\in \cspec H^0(B)$ such that $q\mapsto x$. In this case, $B$ being of finite presentation implies that $\mathbb{L}_B$ has finite \textit{Tor-amplitude}\footnote{We say that a perfect complex $E$ of $ R $-modules has Tor-amplitude in some interval $[a,b]$ if $H^i(E\otimes_R^L N)=0$ for all $i\notin [a,b]$ and for all $R$-modules $N$.}, say in $ [-n,0]$ for some $n\in \mathbb{Z}_+$. Then it follows from \cite[Prop. 7.2.4.23]{Lurie_higheralgebra} that for any integer $k<0$, \ $\mathbb{L}_B[-k-n]$ has Tor-amplitude in $[k,0]$. Notice that both complexes are equivalent. Since we will be interested in $k$-shifted structures, for the  proof, we equivalently use this shifted complex of $B$. Therefore, while getting a refined neighborhood, we assume  w.l.o.g. that for $k<0$,  the  cdga $B$ is such that its cotangent complex $\mathbb{L}_B$ has Tor-amplitude in $[k,0].$

 \end{remark}

\begin{remark}\label{remark_setup before theproof2}
	Let ${\bf X}, x, B, q$ be as in Remark \ref{remark_setup before theproof1}. Then the construction given in \cite[Theorem 4.1]{Brav} (cf. Theorem \ref{localmodelthm}) ensures that there exists a suitable localization of $B$ at $q$, which is equivalent to a minimal standard form cdga $A=A(m)$, for some $m$, constructed inductively as  in (\ref{A(n) construction}) such that there exists $p \in \spec A$ with $p\mapsto q$. Here the integer $m$ is determined by the Tor-amplitude of $\mathbb{L}_B$, which is by assumption $\leq -k$\footnote{It means the perfect complex has Tor-amplitude in $[k,0]$}. It should also be noted that during induction, each $\mathbb{L}_{A(\ell)}$ has Tor-amplitude in $[-\ell,0]$. 
Moreover, as there is an equivalence $A(-k)\rightarrow B$, Proposition \ref{proposition_L as a complex of H^0 modules} provides a simple description for $\mathbb{L}_A$, with Tor-amplitude $\leq -k$. %\footnote{It means the perfect complex has Tor-amplitude in $[k,0]$}. 
%From Proposition \ref{proposition_L as a complex of H^0 modules}, 
That is, when restricted to $ \cspec H^0(A) $,  $\mathbb{L}_A$ is equivalent to the complex $\Omega_A^1 \otimes_A H^0(A)$ of free $H^0(A)$-modules 
\begin{equation*}
	0  \rightarrow  V^k \rightarrow V^{k+1} \rightarrow \cdots \rightarrow V^{-1} \rightarrow V^0 \rightarrow 0,
\end{equation*}with $d|_p^{-i}=0$ for $i=1,2,\dots,-k$ (due to the minimality of $A$). 
\end{remark}
 \pf[\textbf{Proof of Theorem \ref{contact darboux}}]
The proof is essentially  based on that of \cite[Theorem 5.18.]{Brav}.  Let $k<0$ and $x\in \bf X$, apply Theorem \ref{localmodelthm} to get a refined  neighborhood  ${\bf U}= \spec A$ of $x$ with $p \in \cspec H^0(A)$ such that $ i: \spec A \hookrightarrow \bf X $ is an open inclusion, $i(p)=x$, and $A$ is a standard form cdga that is minimal at  $p$. We also assume that $A$ is in fact constructed inductively as described in (\ref{A(n) construction}) with $A=A(-k)$ such that $\mathbb{L}_{A}$ has Tor-amplitude in $[k,0]$. 

  W.l.o.g., we also assume that  $L$ is trivial on $\bf U$, and hence, over $\bf U$, the induced $k$-shifted 1-form $\alpha: \mathbb{T}_X \rightarrow \mathcal{O}_X[k]$ is such that $\mathcal{K}$ is the cocone of ${\alpha}$, up to quasi-isomorphism, and the 2-form $d_{dR}\alpha$ is non-degenerate on $\mathcal{K}$. In that case, the triangle $\mathcal{K} \rightarrow \mathbb{T}_{{\bf X}} \rightarrow L[k]$ splits over $\bf U.$ %We first consider the pullback under the inclusion $i$ of the triangle $\mathcal{K} \rightarrow \mathbb{T}_{{\bf X}} \rightarrow L[k]$, which  splits over ${\bf U}= \spec A$. 
  That is, $\mathrm{cofib}(i^*\mathcal{K} \rightarrow i^*\mathbb{T}_{{\bf X}})\simeq i^*(L[k])$, which is concentrated in $\deg -k$, and $ i^*(\mathbb{T}_{{\bf X}})=  i^*(\mathcal{K})  \oplus i^*(L[k]).$ Here, we will call the 2nd summand $Rest$. Note also that $i^*$ is just the restriction to ${\bf U}=\spec A$, so we will omit it in the notation or use $(-)|_{\spec A}$ instead.
  
  We fix the locally defining 1-form $\alpha$ over $\bf U$ for the rest of the proof.
 %Then the restriction $i^*\alpha$ is a $k$-shifted contact form on $\spec A$.
 % We denote the restriction of $\alpha$ simply by $\alpha_u$. 
 From now on, we will  use the properties of shifted contact structures when restricted to nice local models. %(cf. Remark \ref{remark_cocone/cone vs kernel/cokernel} and relevant discussions after that).
I.e., over $ \spec A $, we use the strict $\ker, \coker $ in $D(Mod_A).$

Consider the sequence $\omega_u:=(d_{dR}\alpha, 0, 0, \dots),$ which defines a closed $k$-shifted 2-form on ${\bf U}$ in the sense of Definition \ref{defn_closed p form}. Applying Proposition \ref{Proposition_exactness} to $k \omega_u$, we obtain  elements $H\in A^{k+1}$ and $\phi \in (\Omega^1_{A})^k$ such that $dH=0$, \ $d_{dR}H+d\phi=0$ , and $k\omega_u \sim (d_{dR}\phi, 0, 0, \dots).$ 
Notice that we in fact have $d_{dR}\phi=kd_{dR}\alpha,$ because there is no non-trivial $\beta \in (\Lambda^2\Omega^1_{A})^{k-1}$ satisfying the relation $ kd_{dR}\alpha- d_{dR}\phi=d \beta$ due to degree reasons.

From Proposition \ref{proposition_L as a complex of H^0 modules}%with $A=A(-k)$
, the tangent complex $\mathbb{T}_{A} |_{\cspec H^0(A)}=  (\mathbb{L}_A|_{\cspec H^0(A)})^{*} $ is also represented by a complex of free finite rank $H^0(A)$-modules, with Tor-amplitude in $[0,-k]$. Then  the $k$-shifted 2-form $d_{dR}\alpha$ defines the induced map of perfect complexes via $v\mapsto \iota_{v}d_{dR}\alpha$:
\begin{equation} 
\begin{tikzpicture}
\matrix (m) [matrix of math nodes,row sep=1em,column sep=1.5em,minimum width=2 em] {
	\mathbb{T}_{A}|_{\cspec H^0(A)}:  & 0  & (V^0)^* & (V^{-1})^* & \cdots & (V^{k+1})^* & (V^k)^* & 0 \\
	\mathbb{L}_A |_{\cspec H^0(A)}[k] :    & 0    &  V^k & V^{k+1} & \cdots & V^{-1} & V^0 & 0,\\
};
\path[-stealth]

%horizontal 1st row
(m-1-2) edge  node [right] { } (m-1-3)
(m-1-3) edge  node [right] { } (m-1-4)
(m-1-4) edge  node [right] { } (m-1-5)
(m-1-5) edge  node [right] { } (m-1-6)
(m-1-6) edge  node [right] { } (m-1-7)
(m-1-7) edge  node [right] { } (m-1-8)
%horizontal 2nd row
(m-2-2) edge  node [right] { } (m-2-3)
(m-2-3) edge  node [right] { } (m-2-4)
(m-2-4) edge  node [right] { } (m-2-5)
(m-2-5) edge  node [right] { } (m-2-6)
(m-2-6) edge  node [right] { } (m-2-7)
(m-2-7) edge  node [right] { } (m-2-8)
%vertical
(m-1-1) edge  node [right] {$ d_{dR}\alpha $} (m-2-1)
(m-1-3) edge  node [right] { } (m-2-3)
(m-1-4) edge  node [right] { } (m-2-4)
(m-1-6) edge  node [right] { } (m-2-6)
(m-1-7) edge  node [right] { } (m-2-7);
%edge [dashed,-] (m-2-1);
\end{tikzpicture}
\end{equation} where both horizontal differentials $d^i, (d^i)^*$ are zero  at $p\in \cspec H^0(A)$ due to the minimality of $A$. 	

 From the contactness condition, $d_{dR}\alpha$ is non-degenerate on the {subcomplex} $\ker \alpha |_{\cspec H^0(A)}$ of $\mathbb{T}_A|_{\cspec H^0(A)}$. %with basis $\{\partial/\partial x_j^{-i}, \partial/\partial y_j^{k+i} \}$. Due to the construction above, this subcomplex can be identified with $\ker \alpha |_{\cspec H^0(A)}$ for some $ \alpha\in (\Omega_A^1)^k $ such that
Therefore, the quotient ${\mathbb{T}_A}/{\ker \alpha}$ can also be identified with a $k$-shifted line bundle $Rest$ (concentrated in degree $-k$) such that one has the  splitting 
\begin{equation} \label{splitting of complexes}
\mathbb{T}_A|_{\cspec H^0(A)}=  \ker \alpha |_{\cspec H^0(A)} \oplus Rest|_{\cspec H^0(A)}.
\end{equation}%where $ \mathbb{T}_A|_{\cspec H^0(A)} $ is the complex of $ H^0(A)$-modules defined by Prop. \ref{proposition_L as a complex of H^0 modules}. %where for each homogeneous degree $i$,   $ (\mathbb{T}_A|_{\cspec H^0(A)})^i=  (\ker \alpha |_{\cspec H^0(A)})^i \oplus Rest^i|_{\cspec H^0(A)}. $

Write $W$ for the \emph{dual} of $ \ker \alpha |_{\cspec H^0(A)}$ in $ \mathbb{L}_A |_{\cspec H^0(A)} $, i.e. $W:= (\ker \alpha |_{\cspec H^0(A)})^*$, then we have the  diagram 

\begin{equation} \label{isomorphism of complexes}
\begin{tikzpicture}
\matrix (m) [matrix of math nodes,row sep=1.5em,column sep=1em,minimum width=1.5 em] {
	W^*\subset \mathbb{T}_{A}|_{\cspec H^0(A)}:  & 0  & (W^0)^* & (W^{-1})^* & \cdots & (W^{k+1})^* & (W^k)^* & 0 \\
	W[k] \subset  \mathbb{L}_A |_{\cspec H^0(A)}[k] :    & 0    &  W^k & W^{k+1} & \cdots & W^{-1} & W^0 & 0\\
};
\path[-stealth]

%horizontal 1st row
(m-1-2) edge  node [right] { } (m-1-3)
(m-1-3) edge  node [right] { } (m-1-4)
(m-1-4) edge  node [right] { } (m-1-5)
(m-1-5) edge  node [right] { } (m-1-6)
(m-1-6) edge  node [right] { } (m-1-7)
(m-1-7) edge  node [right] { } (m-1-8)
%horizontal 2nd row
(m-2-2) edge  node [right] { } (m-2-3)
(m-2-3) edge  node [right] { } (m-2-4)
(m-2-4) edge  node [right] { } (m-2-5)
(m-2-5) edge  node [right] { } (m-2-6)
(m-2-6) edge  node [right] { } (m-2-7)
(m-2-7) edge  node [right] { } (m-2-8)
%vertical
(m-1-1) edge  node [right] {$ d_{dR}\alpha $} (m-2-1)
(m-1-3) edge  node [right] { } (m-2-3)
(m-1-4) edge  node [right] { } (m-2-4)
(m-1-6) edge  node [right] { } (m-2-6)
(m-1-7) edge  node [right] { } (m-2-7);
%edge [dashed,-] (m-2-1);
\end{tikzpicture}
\end{equation} %where, using local coordinates above, \begin{align*}
%W^{-i}&=\langle d_{dR} x^{-i}_1, d_{dR} x^{-i}_2, \dots, d_{dR} x^{-i}_{m_i} \rangle_{A(0)}  \ \text{ for } i=0,1, \dots, \ell.\\
%W^{k+i}&=\langle d_{dR} y^{k+i}_1, d_{dR} y^{k+i}_2, \dots, d_{dR} y^{k+i}_{m_i} \rangle_{A(0)}  \ \text{ for } i=0,1, \dots, \ell, 
%\end{align*} 
such that the vertical maps $ d_{dR}\alpha|_{\ker \alpha} :  (W^{k+i})^* \rightarrow W^{-i}$, $v\mapsto \iota_{v}d_{dR}\alpha$, are  all quasi-isomorphisms. \begin{observation}\label{observation_isomorphisms at p}
	As both horizontal differentials $d^i, (d^i)^*$ are zero  at $p\in \cspec H^0(A)$ because of the minimality of $A$, the vertical maps are isomorphisms at $p$, and hence isomorphisms in a neighborhood of $p$. By localizing $A$ at $p$ if needed, we may assume that the vertical maps are all \emph{isomorphisms.} 
\end{observation} 

\paragraph{When $k$ is odd.} We now focus on a particular and the simplest case: $k$ is \textit{odd}. %As noted before, similar local models can be obtained for the other cases $ (a) \ k\equiv0 \mod 4, \text{ and } (b) \ k\equiv2 \mod 4$ by modifying graded variables, the Hamiltonian, and the internal differential. %Therefore, depending on the given $k<0$, the same arguments below will also work after replacing the graded variables and the expressions for $H, \phi, d$ by the suitable analogous ones that are explicitly given in Equations (\ref{new local variables for k=-4l}) $ - $ (\ref{defn_phiv2}).
Let $k=-2\ell -1$ for  $\ell\in\N$. Localizing $A$ at $p$ if necessary, first choose 
degree 0 variables $x_1^0, x_2^0, \dots, x_{m_0}^0$ in $A(0)$ such that $\{d_{dR}x_j^0: j=1,\dots,m_0 \}$ forms a basis for  $W^0$ over $A(0),$ and   $(Rest^*)^0=0$

%Next, for $i=1,\dots, \ell$, choose $ x_1^{-i}, x_2^{-i}, \dots, x_{m_i}^{-i} \in A^{-i}$ such that  $d_{dR}x_1^{-i}, \dots d_{dR}x_{mi}^{-i}$ form a basis of $W^{-i}$ over $A(0)$, and  $(Rest^*)^{-i}$ is trivial over $A(0)$. 

Now, by the equivalence $ d_{dR}\alpha :  (W^{k+i})^* \rightarrow W^{-i}$, we have $H^{-i}(W)\simeq H^{k+i}(W)^*,$ and hence $\dim H^{-i}(W) = \dim H^{k+i}(W).$ Then we can choose the following set of generators for the free algebra $A$:% $A$ is free over $A(0)$ with $m_i$ generators in degree $-i$ for $i=1, \dots, \ell$, and $m_i$ generators in degree $k+i$ for $i=1, \dots, \ell$. 
\begin{itemize}
	\item When $i=0$, we can find generators $ y_1^{k}, y_2^{k}, \dots, y_{m_{0}}^{k}, z^{k}  \in A^{k}$ such that $ \{ d_{dR}y_1^{k}, \dots, d_{dR}y_{m_0}^{k}\}$ is a basis for $W^{k}$  which is dual to  the basis $\{ d_{dR} x_1^{0}, \dots ,d_{dR} x_{m_0}^{0}\}$ for $W^{0}$ and that the complex $Rest$, which is concentrated in $\deg -k$, is  generated by the vector field $\partial/\partial z^{k}$ of degree $-k$.  
	\item Similarly, for $1 \leq i\leq \ell$, choose  $x_1^{-i}, \dots, x^{-i}_{m_i} \in A^{-i}; \text{ and } y_1^{k+i}, \dots, y_{m_i}^{k+i}  \in A^{k+i}$ such that $ \{ d_{dR}y_j^{k+i} : j= 1,2, \cdots, m_i \}$ is a basis for $W^{k+i}$ over $A(0)$ which is dual to  the basis $\{ d_{dR}x_1^{-i}, \dots ,d_{dR}x_{m_i}^{-i}\}$ for $W^{-i}$.
\end{itemize} That is,  using the local coordinates above,  $ \text{for } i=0,1, \dots, \ell, $ we have \begin{align*}
W^{-i}&=\langle d_{dR} x^{-i}_1, d_{dR} x^{-i}_2, \dots, d_{dR} x^{-i}_{m_i} \rangle_{H^0(A)},  \\
W^{k+i}&=\langle d_{dR} y^{k+i}_1, d_{dR} y^{k+i}_2, \dots, d_{dR} y^{k+i}_{m_i} \rangle_{H^0(A)}.
\end{align*} By Observation \ref{observation_isomorphisms at p},  the isomorphisms in Diagram (\ref{isomorphism of complexes}) imply that for $ i=0,1, \dots, \ell, $ we have
\begin{align}
(W^{-i})^*&=\langle \partial/\partial x_1^{-i}, \dots, \partial/\partial x_{m_i}^{-i} \rangle_{H^0(A)}  \xrightarrow{\sim} \langle d_{dR} y^{k+i}_1, \dots, d_{dR} y^{k+i}_{m_i} \rangle_{H^0(A)}, \label{sending duals 1}\\
(W^{k+i})^*&=\langle \partial/\partial y^{k+i}_1, \dots, \partial/\partial y^{k+i}_{m_i} \rangle_{H^0(A)}  \xrightarrow{\sim} \langle d_{dR} x_1^{-i}, \dots, d_{dR}x^{-i}_{m_i} \rangle_{H^0(A)}. \label{sending duals 2} 
\end{align} Then, by the splitting in Equation (\ref{splitting of complexes}), we have
\begin{align} \label{explicit generators for ker and rest}
\ker \alpha |_{\cspec H^0(A)} &=\big\langle \partial/\partial x_1^{-i}, \dots, \partial/\partial x_{m_i}^{-i},  \partial/\partial y^{k+i}_1, \dots, \partial/\partial y^{k+i}_{m_i} : i=0,1, \dots, \ell \ \big\rangle_{H^0(A)},  \nonumber \\
Rest|_{\cspec H^0(A)} &=\big\langle \partial/\partial z^k \big\rangle_{H^0(A)} . 
\end{align} Here $ Rest|_{\cspec H^0(A)} $ is a subcomplex of $\mathbb{T}_A|_{\cspec H^0(A)}$ that is concentrated in degree $-k$ and  generated by the degree $-k$ vector field $\partial/\partial z^{k}$. Moreover, $ \ker \alpha$ has $\mathrm{Tor}$-amplitude $[0,-k]$.

Using the splitting (\ref{splitting of complexes}) with the explicit forms (\ref{explicit generators for ker and rest}), the non-degeneracy condition on $d_{dR}\alpha |_{\ker \alpha}$  sending dual basis of $ d_{dR} x_1^{-i}, \dots, d_{dR}x^{-i}_{m_i} $ to the basis $  d_{dR} y^{k+i}_1, \dots, d_{dR} y^{k+i}_{m_i} $ (and vice versa) as in Equations (\ref{sending duals 1})- (\ref{sending duals 2}) implies that $\dR \alpha \in \wedge^2\Omega_A^1[k]$ must be of the form
\begin{equation}\label{2form coming from non-degeneracy}
d_{dR}\alpha =\displaystyle \sum_{i=0}^{\ell} \sum_{j=1}^{m_i} d_{dR}x_j^{-i} d_{dR}y_j^{k+i}.
\end{equation} 
Notice  that the kernel $\ker (d_{dR}\alpha \cdot)$ of the induced morphism $d_{dR}\alpha \cdot: \mathbb{T}_A \rightarrow \Omega^1_A[k]$ is spanned by the degree $-k$ vector field $\partial/\partial z^k$, while $ d_{dR}\alpha \cdot $ acts on $\ker \alpha$ as  in Equations (\ref{sending duals 1})-(\ref{sending duals 2}). So,  we can get the identification $\faktor{\mathbb{T}_{{A}}}{\ker d_{dR}\alpha}\simeq  \ker \alpha$.

Scaling $z^k$ we may assume $\iota_{\partial/\partial z^k} \alpha =1$. Now, our goal is to find a unique  $\alpha$ satisfying Eqn. (\ref{2form coming from non-degeneracy}), the condition on the kernel in (\ref{explicit generators for ker and rest}), and the equation $\iota_{\partial/\partial z^k} \alpha =1 $. %Since the  form $ \alpha $ can also been seen as a representative of $\phi$, we may let $ \alpha=\dR z^k + \phi/k $. 
Then  such  $ \alpha $ satisfying the desired properties (uniqueness discussed later) can be explicitly written as  \begin{align}\label{desired form}
\alpha
&= d_{dR}z^k+  \sum_{i,j} y_j^{k+i}d_{dR}x_j^{-i}.
\end{align} %such that $\dR \alpha$ recovers (\ref{2form coming from non-degeneracy}). Note also that the kernel of $\alpha$ for the above representation is generated by the vectors fields\footnote{The vector fields $ ix^{-i}_j\partial/\partial x_j^{-i} + (k+i)y^{k+i}_j\partial/\partial y^{k+i}_j $ also annihilate $\alpha$, but they can be expressed using the others.} of the form  $ k \partial/\partial x_j^{-i}-(k+i)y^{k+i}_j\partial/\partial z^{k} \text{ and } k\partial/\partial y^{k+i}_j+ix^{-i}_j\partial/\partial z^{k} $%$ \partial/\partial x_j^{-i} \text{ and } \partial/\partial y^{k+i}_{j}-x^{-i}_{j} \partial/\partial z^k $
%, which can be identified with the ones in (\ref{explicit generators for ker and rest}).   

\vspace{2pt}

\bfem{To sum up,} $A$ is  identified as a \emph{commutative graded algebra} with the commutative graded algebra over $A(0)$ freely generated by the variables $z^k, x_j^{-i}, y_j^{k+i}$ exactly as in Example \ref{model example}, together with the representation of $\alpha$   above.
It remains to show that  the differential $d$ can be given by Eqn. (\ref{defn_internal d contact}). To this end, we analyze the defining equations for  the pair $(H,\phi)$.

First of all, Brav, Bussi, and Joyce \cite[Theorem 5.18]{Brav}  show in their proof that combining  the  equation $d_{dR}\phi=kd_{dR}\alpha$ with the formula (\ref{2form  coming from non-degeneracy}), we may explicitly write\footnote{Alternatively,  one can let $ \phi=  k\sum_{i=0}^{\ell} \sum_{j=1}^{m_i} y_j^{k+i}d_{dR}x_j^{-i}$.   Leaving $\dR\phi$ unchanged, these expressions can be transformed to each other by replacing $H, \phi$ by suitable $H+d(\cdots), \ \phi+d_{dR}(\cdots),$ respectively.} $$ \phi= \sum_{i=0}^{\ell} \sum_{j=1}^{m_i} [(-i)x_j^{-i}d_{dR}y_j^{k+i}+(k+i)y_j^{k+i}d_{dR}x_j^{-i}].$$%we start with replacing $H,\phi$ by \begin{align*}
	%H'&=H-d\big[ \sum_{i,j} (-1)^{i} ix_j^{-i} y_j^{k+i} \big] \\ 
	% \phi'&= \phi-\dR\big[ \sum_{i,j} (-1)^{i} ix_j^{-i} y_j^{k+i} \big] = \sum_{i,j} \big[(-i)x_j^{-i}d_{dR}y_j^{k+i}+(k+i)y_j^{k+i}d_{dR}x_j^{-i}\big].
%\end{align*} Note that $H',\phi'$ also satisfy the defining equations. Then 
%modifying\footnote{Use $H-d\big[ \sum_{i,j} (-1)^{i} ix_j^{-i} y_j^{k+i} \big] $ and $\phi-\dR\big[ \sum_{i,j} (-1)^{i} ix_j^{-i} y_j^{k+i} \big]$, which gives the form in Footnote 7.} $H,\phi$ properly, 
Then \cite[Theorem 5.18]{Brav} also shows that expanding %\footnote{While expanding, write $\phi= \sum_{i,j} [(-i)x_j^{-i}d_{dR}y_j^{k+i}+(k+i)y_j^{k+i}d_{dR}x_j^{-i}] + \dR [ \sum_{i,j} (-1)^{i} ix_j^{-i} y_j^{k+i} ]$ and replace $H$ by $H+d\big[ \sum_{i,j} (-1)^{i} ix_j^{-i} y_j^{k+i} \big]. $} 
$ d_{dR}H+d\phi=0 $ with the explicit representation of $\phi$ above and comparing the coefficients of $d_{dR}$-terms, one gets the following formulas for $d$:
\begin{equation} 
	d|_{A(0)}=0; \ dx_j^{-i} =  \dfrac{\partial H}{\partial y_j^{k+i}} \text{ for all } i>0,j; \ \text{ and } \ dy_j^{k+i} =  \dfrac{\partial H}{\partial x_j^{-i}} \text{ for all } i,j.
\end{equation}  

Next, using these equations for the differential to expand\footnote{These equations will be enough as $H$ is independent of the top degree variables $z^k, y^k_j$'s.} $dH=0$, \cite[Theorem 5.18]{Brav} also implies that $H$ must satisfy the classical master equation (\ref{defn_CME}). 
Before the final step, we also observe that using the explicit representation of $\phi$ above, $\alpha$ in (\ref{desired form}) can also be rewritten as \begin{equation} \label{form written alternatively}
	\alpha=\dR z^k + \dfrac{1}{k}\left[\phi + \dR  \left(\sum_{i,j} (-1)^{i} ix_j^{-i} y_j^{k+i} \right)\right].
\end{equation} 

Finally, combining the equation $d_{dR}H+d\phi=0 $ with the version of $\alpha$ in (\ref{form written alternatively}) (and $d\alpha=0$), we get $ \dR H=-d\phi=-kd\alpha + kd\circ\dR z^k +d\circ\dR [\cdots]= \dR (-kdz^k - d[\cdots]).  $ So, we  have  such $ H $ satisfying $ -kdz^k=H+d[\cdots]$. It follows  that the differential $d$ is fully given as in  (\ref{defn_internal d contact}) and that $(A,d)$ is identified  with the cdga in Example \ref{model example}.

We then conclude that $(A,d,\alpha)$ is in contact Darboux form, which completes the proof of Theorem \ref{contact darboux} when $k<0$ is odd.

\paragraph{Uniqueness.} Assume that there is an element $\alpha'$ satisfying the desired properties listed above. Let us clarify  how the corresponding conditions uniquely (up to interchange %interplay/choice of roles
of $x_{j}^{-i} \text{ and } y_{j}^{k+i}$)\footnote{The roles of $x_{j}^{-i},y_{j}^{k+i}$ are symmetric in (\ref{2form coming from non-degeneracy}) and (\ref{explicit generators for ker and rest}), where $ \dR x_{j}^{-i}\dR y_{j}^{k+i}= \dR y_{j}^{k+i}\dR x_{j}^{-i} $ for $k$ odd.} determine the  representation in (\ref{desired form}). %Using  Eqn. (\ref{2form coming from non-degeneracy}) and the coordinates $x^{-i}_j, y^{k+i}_j, z^k$, we can write $$ \alpha'= \sum_{i,j} \left[a^{k+i}_j \dR x^{-i}_j + b^{-i}_j \dR y^{k+i}_j + c^0\dR z^k\right],$$ where   $a_{\nu}^{\mu},b_{\nu}^{\mu} \in A^{\mu}, \  c^0\in A^0$ such that  $b_{\nu}^{\mu}$'s depend only on $x_{j'}^{-i'}$'s and  $A(0)$ for degree reasons. %we have $ kd_{dR}\alpha' = k\sum_{i,j}d_{dR}x_j^{-i} d_{dR}y_j^{k+i}$. %and $\ker (d_{dR}\alpha' \cdot)= \langle \partial/\partial z^k\rangle,$ such that ${\mathbb{T}_{{A}}}/{\ker d_{dR}\alpha'}\simeq  \ker \alpha'$.  
We first observe that due to Eqn. (\ref{2form coming from non-degeneracy}) and the  condition  $\iota_{\partial/\partial z^k} \alpha' =1,$ any such element
$\alpha'$ takes the form %$\dR z^k +  \psi$ for some $\psi \in \Omega_A^1[k]$ given by 
\[ \alpha'= \dR z^k +  \psi \quad \text{with} \quad \psi = \sum_{i,j} \left[a^{k+i}_j \dR x^{-i}_j + b^{-i}_j \dR y^{k+i}_j\right], \] where   $a_{\nu'}^{\mu'}\in A^{\mu'}, b_{\nu}^{\mu} \in A^{\mu}$ such that  $b_{\nu}^{\mu}$'s depend only on $x_{j'}^{-i'}$'s and  $A(0)$ for degree reasons.
From the condition  (\ref{explicit generators for ker and rest}), we take\footnote{Instead, one may set $ a^{k+i}_j=0 \ \forall i,j$, giving rise to the form of $\alpha'$ with the roles of $ x_{\nu}^{\mu},y_{\nu}^{\sigma} $ interchanged.  %Either version of $k\alpha'-k\dR z^k$ induces alternative solutions to the defining equations, see Remark \ref{remark_alternative and general versions of alpha}. Note also that by \cite[Prop. 5.7(c)]{Brav}, they differ by a $d[\cdots]$ term, meaning that they are equivalent, see Defn \ref{defn_p form of deg k}. 
} 
$ b^{-i}_j=0$ for all $i,j$, %(which also amounts to replacing $\psi$ by a suitable $\psi-\dR [\cdots]$)
leading to 
\[ d_{dR}\alpha' = \dR \psi=\sum_{i,j}d_{dR}a_j^{k+i} d_{dR}x_j^{-i}.\] 
Now, using the RHS of (\ref{2form coming from non-degeneracy}) for comparison, %(see also \cite[Proof of Thm. 5.18]{Brav}), observe  that $ a_j^{k+i}= y_j^{k+i} + (\deg \geq2 \text{ terms in }x_{j'}^{-i'} \text{ for } i>0 \text{ and } y_{j''}^{k+i''})$.Hence, $ a_j^{k+i}$ are alternative choices for $y_j^{k+i}$. 
we have %\footnote{Such element is in fact of the form $ a_j^{k+i}= y_{j}^{k+i} + (\deg \geq2 \text{ terms in }x_{j'}^{-i'} \text{ for } i>0 \text{ and } y_{j''}^{k+i''})$.} 
$ a_j^{k+i}= y_j^{k+i}$ for all $i,j$. Replacing them accordingly, we obtain $ \psi = \sum_{i,j} y^{k+i}_j \dR x^{-i}_j,$ which gives $\alpha'= \dR z^k +  \sum_{i,j} y^{k+i}_j \dR x^{-i}_j,$ hence the desired form (\ref{desired form}).

%Fixing $(H,\phi)$ above, assume that there is $\alpha'$ satisfying the desired properties. Since $\dR \phi = k\dR \alpha'$, we can write $k\alpha'=\phi + \dR \theta$ for some $\theta\in A^k$. Due to the  conditions (\ref{explicit generators for ker and rest}) and  $\iota_{\partial/\partial z^k} \alpha' =1,$ the term  $ \dR \theta $ takes the form $k\dR z^k + \dR \theta'$, where $\theta'$ is a degree $k$ element not involving $z^k.$ Thus, we write $k\alpha'=k\dR z^k +\phi + \dR \theta'$. Now, add suitable $\dR [\cdots]$ on both sides to get $k\alpha'+\dR [\cdots]=k\alpha+ \dR \theta'$, meaning that $ \alpha-\alpha'$ is $\dR$-exact. 

\paragraph{When $k$ is not odd.} For the other cases $ (a) \ k\equiv0 \mod 4, \text{ and } (b) \ k\equiv2 \mod 4  $, one should use another sets of variables as in Equations (\ref{new local variables for k=-4l}) and (\ref{new local variables for k=-4l-2}), respectively, and modify $H, \phi, d$ as in Equations (\ref{new local variables for k=-4l}) $ - $ (\ref{defn_phiv2}).  We leave details to the reader.

This completes the proof of Theorem \ref{contact darboux} for all $k<0$, and hence that of Theorem \ref{THM1}.

\epf

\begin{observation}
	Denote by $B^0$ the subalgebra of $A^0$ with basis $ x_1^0, x_2^0, \dots, x_{m_0}^0, $ Then we define a sub-cdga $B$ of $A$ to be the free algebra over $B^0$ on generators $x_j^{-i}, y_j^{k+i}$ only, with inclusion $ \iota: B \hookrightarrow A. $  Observe that   the elements $\phi, \omega^0:=d_{dR}\alpha |_{\ker \alpha}$ are all images under $\iota$ of the elements $\phi_B, \omega_B^0:=\sum_{i=0}^{\ell} \sum_{j=1}^{m_i} d_{dR}x_j^{-i} d_{dR}y_j^{k+i}$, respectively. %Moreover, we have $\spec A=U\xrightarrow{j:=\spec (\iota)} V:=\spec B$ such that $j^*\omega_B^0 \sim d_{dR}\alpha_u |_{\ker \alpha_u}$. Then $ \omega_B:= (\omega^0_B,0,0, \cdots) $ is a $k$-shifted symplectic structure on $\spec B$.  
	As in Section \ref{the pair},  $B$ is a minimal standard form cdga which in fact serves as a local symplectic  model (for $k<0$ odd).  As noted before, similar local models can be explicitly obtained for the other cases using Equations (\ref{new local variables for k=-4l}) $ - $ (\ref{defn_phiv2}).  
	
	Suppose that we construct such  $(A,B)$ for $k<0$ with $k \not\equiv 2 \mod 4 $. From Observation \ref{observation_vdim symplectic},  the virtual dimension $\vdim B$ is then always \textit{even}, and hence  $\vdim A= \vdim B \pm1$ is \emph{odd}. Otherwise, $\vdim A$ can take any value in $ \mathbb{Z} $. That is, if  $A$ is a cdga in contact Darboux form, we have 
	 \[\vdim A = 
	\begin{cases*} 
		-1, & \text{ if }$ k $ \mbox{is odd},\\
		\text{odd in } \mathbb{Z}, & \text{ if }$ k/2 $ \mbox{is even}, \\
		\text{any value in } \mathbb{Z}, & \text{ if } $ k/2 $ \mbox{is odd}. \\
	\end{cases*}
	\]
\end{observation}

\newpage
\section{Symplectifications} \label{section_symplectization}

In this section, we give the formal description of  the \textit{symplectification} of a negatively shifted contact derived $\K$-scheme.

\subsection{Symplectified spaces for contact manifolds}In classical contact geometry,  for a contact manifold $(M,\xi= \ker (\alpha))$ with a globally defined contact 1-form $\alpha$, one can define the \textit{symplectification} $\widetilde{M}$ of $M$ as the total space of the  bundle $M \times \mathbb{R}^*\rightarrow M$ with a canonical symplectic form $\omega:=d_{dR}(e^t \alpha)$, where $t$ is  the $\mathbb{R}$-coordinate. Most of the standard references \cite{ACSilva,Geiges} use this approach, where the contact structure $\xi$ is in fact assumed to be coorientable (see Remark \ref{remark_coorientible and global alpha} for the definition). 

However, for non-coorientible contact structures, the coordinate-dependent description above can no longer be applicable (as no global $\alpha$ and $t$-variable available); instead, we may use the following description from \cite[Appendix 4.E]{Arnold}: \ Given a contact manifold $(M, \xi)$, we let
\begin{equation} \label{defn_classical symplectizastin Arnold style}
\widetilde{M}:= \big\{(p, \alpha_p ) :  p \in M, \ \alpha_p\in T_p^*M, \ \text{s.t. }\ker \alpha_p = \xi_p\big\}.
\end{equation}Here, $\widetilde{M}$ is just the set of all contact forms on the contact manifold. It should be noted that for a pair $(p, \alpha_p)\in \widetilde{M}$, $ \ \alpha_p $ is not a differential form but just a linear form on one tangent space $T_pM$ at the point of contact of the manifold such that its zero set is the contact plane. From \cite[Appendix 4.E]{Arnold}, it is straightforward to see that $ \widetilde{M} $ is a smooth manifold of even dimension  $\dim M +1$. Notice that there is a natural $\mathbb{R}^*$-action on $\widetilde{M}$ via $f\cdot (p, \alpha_p )= (p, f\alpha_p )$ such that $ \widetilde{M}/\mathbb{R}^*\simeq M. $ Therefore, $\widetilde{M}$ can be identified as the total space of the $\mathbb{R}^*$-bundle over $M$.

From this identification with the $\mathbb{R}^*$-bundle $\pi: \widetilde{M}\rightarrow M$, the canonical \textit{symplectic structure} on $ \widetilde{M} $ can be defined as $\omega:= d_{dR}\lambda$, where  the so-called \textit{canonical 1-form} $\lambda$ is the differential 1-form on $\widetilde{M}$ whose value on any vector $v\in T_x\widetilde{M}$ at a point $x=(p,\alpha_p) \in \widetilde{M}$ is given by \begin{equation} \label{defn_canonocal 1form_classical}
\lambda_x(v):= \alpha_p(\pi_{*,x}(v)).
\end{equation}

\subsection{Symplectified derived spaces}The construction of a canonical symplectified space associated to a contact space in terms of the total space of a line bundle with $\R^*$-action can be promoted to derived symplectic geometry. This leads to the definition of the \emph{symplectification} of a shifted contact derived stack. In what follows, we explain the details. Let us begin by some relevant notions.
	
	\begin{definition}
		Let ${\bf X}\in dStk_{\K}$ and $E\in QCoh(\bf X) $, then the \emph{total space $\bf \widetilde{E}$ of $E$} is defined as a derived stack sending  \begin{equation}
		A \mapsto \widetilde{{\bf E}}(A):=\{(p,s) : \ p\in {\bf X}(A), \ s\in p^*E    \}, \text{ with }  A\in cdga_{\K}.
		\end{equation} Here, from Yoneda's lemma, ${\bf X}(A) \simeq Map_{dPstk}(\spec A, {\bf X})$, and hence any $A$-point $p\in {\bf X}(A) $ can be seen as a morphism $p: \spec A \rightarrow \bf X$ of derived pre-stacks, and thus its pullback map $ p^*: QCoh({\bf X}) \rightarrow QCoh(\spec A) \simeq Mod_A$ sends $E \mapsto p^*E$. Hence, $s$ is just an element of the $A$-module $ p^*E,$ a "fiber" over $p$.
	\end{definition}

\begin{example}
 if ${\bf X}\in dStk_{\K}$ admits a cotangent complex (which is always the case when $\bf X$ is also Artin), we can define the cotangent stack $\bf T^*X$ to be the total space $\bf \widetilde{{L}_{{\bf X}}}$ of $\mathbb{L}_{{\bf X}}\in QCoh(\bf X)$ and the $n$-shifted cotangent stack $\bf T^*[k]X$ to be the total space $\widetilde{\bf{L}_{{\bf X}}[k]}$ of $\mathbb{L}_{{\bf X}}[k]\in QCoh(\bf X)$. For more details see \cite[\S 2]{Safronov}. 
\end{example}

Recall from \cite{Damien2} that for a perfect module $E$ over $\bf X$, its stack of sections $ \widetilde{{\bf E}}$, defined by $\widetilde{{\bf E}}(-)=\mathbb{R}\spec_{{\bf X}}(Sym(E^{\vee}))(-)$ is acted on by $\mathbb{G}_m$ because $Sym(E^{\vee})$ is graded $\mathcal{O}_{{\bf X}}$-algebra. This new grading is then called the \emph{fibre grading}. Note also that the both zero section $\bf X \rightarrow \widetilde{E}$ and the projection $\bf \widetilde{E} \rightarrow X$ are  $\mathbb{G}_m$-equivariant for the trivial $\mathbb{G}_m$-action on $\bf X.$ With this terminology, $\bf T^*[k]X$ is nothing but the stack of $ k$-shifted 1-forms on $ \bf X$ with a natural $\mathbb{G}_m$-action. 

 More generally, for $E\in QCoh(\bf X) $, the $\mathbb{G}_m$-action is given as a morphism of derived $\K$-stacks $\triangleleft: \mathbb{G}_m \times_{\cspec \K} \bf \widetilde{E} \rightarrow \widetilde{E}$ such that for each $A\in cdga_{\K}$, we have a $\mathbb{G}_m(A)$-action $\triangleleft_A$ on ${\bf \widetilde{E}}(A)$,
where $\mathbb{G}_m$ is the functor that maps $A \mapsto A^{\times}$. This definition also holds for any derived $S$-stack.

Now, we are in place of introducing the definition of the \emph{symplectification} of a $k$-shifted contact derived  stack using the machinery above:

\begin{definition} \label{defn_symplectification}
Let $\bf X$ be a locally finitely presented derived $\K$-scheme with a $k$-shifted contact structure $(\mathcal{K}, \kappa, L)$. The \bfem{symplectification} is the total space $\widetilde{\bf L}$ of the $\mathbb{G}_m$-bundle of $L$, with a canonical $k$-shifted symplectic structure (for which the $\mathbb{G}_m$-action is of weight 1) as defined below. 
\end{definition}

\subsubsection{Step-1:  The derived stack $ \mathcal{S}_{{\bf{X}}} $ of contact forms} Let $({\bf X}; \mathcal{K}, \kappa, L )$ be a $k$-shifted contact derived $\K$-scheme of locally finite presentation. Given $k<0$ and $p\in \bf X$, find a minimal standard form cdga $ A $ and an affine derived sub-scheme ${\bf U}:= \spec A$ such that $p: \spec A \rightarrow \bf X$ is Zariski open inclusion. %(we may further assume $A$ is of minimal standard form)
 Here, we assume w.l.o.g. that  $L$ is trivial on $\bf U$. %Choose a $k$-shifted contact form on $\bf U$  such that $Cone(\alpha)\simeq \mathcal{K}^{\vee}[k]$. 

Define a functor $ \mathcal{S}_{{\bf{X}}}: cdga_{\K} \rightarrow Spcs$ by $ A\mapsto \mathcal{S}_{{\bf{X}}}(A), $ where
\begin{equation}\label{defn_presstack S_X(A)}
 \mathcal{S}_{{\bf{X}}}(A):= \big\{(p, \alpha, v ) :  p \in {\bf X} (A), \ \alpha:   p^*(\mathbb{T}_{{\bf X}})\rightarrow \mathcal{O}[k], \ v: Cocone(\underline{\alpha}) \xrightarrow{\sim}p^*(\mathcal{K}) \big\},
\end{equation} where each $ v $ is a quasi-isomorphism respecting the natural morphisms $ p^*\kappa: p^*\mathcal{K} \rightarrow p^*(\mathbb{T}_{{\bf X}})$ and $ Cocone(\underline{\alpha}) \rightarrow p^*(\mathbb{T}_{{\bf X}})   $.
Under the current assumptions, the perfect complexes $\mathbb{T}_A, \mathbb{L}_A$, when restricted to ${\cspec H^0(A)}$, are both quasi-isomorphic to free  complexes of $H^0(A)$-modules. For $A \in cdga_{\K}$, we then define a $\mathbb{G}_m(A)$-action on $ \mathcal{S}_{{\bf{X}}}(A)  $ by 
\begin{equation*}
f \triangleleft (p, \alpha, v):= (p, f\cdot \alpha,v).
\end{equation*} %Consider the functor sending $A \mapsto H^0(A).$ Denote the image under $H^0$ of an element $f$ simply by $f^0.$ Localizing $ A $, if necessary, w.l.o.g. we may assume that the image $f^0$ is always invertible. It follows that $f^0$ lies in $(A^0)^{\times}$, which is by definition $\mathbb{G}_m(A^0)=(A^0)^{\times}$.

Now, the following observation endows $ \mathcal{S}_{{\bf{X}}}$ with the structure of a derived stack:

\begin{proposition}
	$ \mathcal{S}_{{\bf{X}}}$ is equivalent to the total space $\bf \widetilde{L}$ of the $\mathbb{G}_m$-bundle of $L$\footnote{In general, when the group scheme $G=GL_n$, there is an equivalence between locally free sheaves of  $\rank n$ and $GL_n$-torsors (hence principal $GL_n$-bundles). 
		In this regard, a line bundle $L$ is nothing but a $\mathbb{G}_m$-bundle.}. Therefore, it has the structure of a derived stack together with the projection map $\pi:  \mathcal{S}_{{\bf{X}}}\rightarrow \bf X$.  We then call $ \mathcal{S}_{{\bf{X}}} $ the \bfem{derived stack of $k$-contact forms}.
\end{proposition}
\pf
Let $ (p, \alpha, v ) $ be a point in $ \mathcal{S}_{{\bf{X}}}(A)$. From definitions, we have two homotopy fiber sequences $(i) \ Cocone(\underline{\alpha})\rightarrow p^*(\mathbb{T}_{{\bf X}})  \rightarrow \Ima \alpha$ and $(ii) \ Cocone({\alpha})\rightarrow p^*(\mathbb{T}_{{\bf X}}) \rightarrow \mathcal{O}[k]$. Then we get the following observation:

\begin{observation} \label{new triangle}
There exists a triangle $ Cocone(\underline{\alpha})\rightarrow p^*(\mathbb{T}_{{\bf X}})  \rightarrow \mathcal{O}[k]$.
\end{observation}

\pf[Proof of Observation \ref{new triangle}] From the first sequence $ (i) $ above, we have  the following homotopy commutative diagram, where the left-hand square is the homotopy fiber. 
\begin{equation} 
	\begin{tikzpicture}
		\matrix (m) [matrix of math nodes,row sep=2em,column sep=2em,minimum width=2 em] {
			Cocone(\underline{\alpha})	& \star  &  \star  \\			
			p^*(\mathbb{T}_{{\bf X}})  & \Ima \alpha    & \mathcal{O}[k] \\
		};
		\path[-stealth]
		
		%horizontal 1st row
		(m-1-1) edge  node [above] { $ \pi_2 $ } (m-1-2)
		(m-1-2) edge  node [above] { {\small id} } (m-1-3)
		%(m-1-4) edge  node [right] { } (m-1-5)
		%(m-1-5) edge  node [right] { } (m-1-6)
		%(m-1-6) edge  node [right] { } (m-1-7)
		%(m-1-7) edge  node [right] { } (m-1-8)
		%horizontal 2nd row
		(m-2-1) edge  node [above] { $ \underline{\alpha} $ } (m-2-2)
		(m-2-2) edge  node [above] { $ j $ } (m-2-3)
		%(m-2-4) edge  node [right] { } (m-2-5)
		%(m-2-5) edge  node [right] { } (m-2-6)
		%(m-2-6) edge  node [right] { } (m-2-7)
		%	(m-2-7) edge  node [right] { } (m-2-8)
		%vertical
		(m-1-1) edge  node [left] {$ \pi_1 $ } (m-2-1)
		(m-1-2) edge  node [right] { $ 0 $ } (m-2-2)
		(m-1-3) edge  node [right] {$ 0 $ } (m-2-3)
		%(m-1-6) edge  node [right] { } (m-2-6)
		%(m-1-7) edge  node [right] { } (m-2-7);
		%edge [dashed,-] (m-2-1)
		;
	\end{tikzpicture} 
\end{equation} There is a homotopy $H$ with $\underline{\alpha}\circ \pi_1 \sim 0$. Using the monomorphism $j$, we get a homotopy $j\circ H$ such that $  j\circ\underline{\alpha}\circ \pi_1 \sim 0$, where $j\circ\underline{\alpha}\sim \alpha$. It follows that $ Cocone(\underline{\alpha}) $ homotopy commutes the outer diagram. Then the universality of $ Cocone({\alpha}) $ implies that  $\exists \ \varphi_1: Cocone(\underline{\alpha})\rightarrow Cocone({\alpha}) $.

Likewise, from the second sequence $ (ii) $ above, we have  the following homotopy commutative diagram. 
\begin{equation} 
	\begin{tikzpicture}
		\matrix (m) [matrix of math nodes,row sep=1.8em,column sep=2em,minimum width=2 em] {
			Cocone({\alpha})	&   &  \star  \\	
								& \star & \\				
			p^*(\mathbb{T}_{{\bf X}})  &     & \mathcal{O}[k] \\						
		                               & \Ima \alpha    &  \\
		};
		\path[-stealth]
		
		%horizontal 1st row
		(m-1-1) edge  node [above] { $ \pi_2 $ } (m-1-3)
		(m-1-1) edge  node [right] { $ \pi_2 $ } (m-2-2)
		(m-2-2) edge  node [below] { \small{id} } (m-1-3)
		%(m-1-4) edge  node [right] { } (m-1-5)
		%(m-1-5) edge  node [right] { } (m-1-6)
		%(m-1-6) edge  node [right] { } (m-1-7)
		%(m-1-7) edge  node [right] { } (m-1-8)
		%horizontal 2nd row
		(m-2-1) edge  node [above] { $ \alpha $ } (m-2-3)
		%(m-2-2) edge  node [above] { $ i $ } (m-2-3)
		%(m-2-4) edge  node [right] { } (m-2-5)
		%(m-2-5) edge  node [right] { } (m-2-6)
		%(m-2-6) edge  node [right] { } (m-2-7)
		%	(m-2-7) edge  node [right] { } (m-2-8)
		%vertical
		(m-1-1) edge  node [left] {$ \pi_1 $ } (m-2-1)
		(m-2-2) edge  node [right] { $ 0 $ } (m-4-2)
		(m-1-3) edge  node [right] {$ 0 $ } (m-2-3)
		(m-3-1) edge  node [below] {$ \underline{\alpha} $ } (m-4-2)
		(m-4-2) edge  node [below] {$ j $ } (m-3-3)
		%(m-1-6) edge  node [right] { } (m-2-6)
		%(m-1-7) edge  node [right] { } (m-2-7);
		%edge [dashed,-] (m-2-1)
		;
	\end{tikzpicture} 
\end{equation} Here, we have a homotopy such that $ \alpha \circ \pi_1 \sim 0. $ From the epi-mono factorization of $\alpha$, we have $  j\circ\underline{\alpha}\circ \pi_1 \sim 0$ as well. From $j\circ 0 \sim 0$, we obtain $j\circ\underline{\alpha}\circ \pi_1 \sim j \circ 0.$ Since $j$ is a monomorphism, we get the induced homotopy such that $\underline{\alpha}\circ \pi_1 \sim  0.$ It follows from the universality of $ Cocone(\underline{\alpha}) $ that there exists a map $ \varphi_2: Cocone({\alpha}) \rightarrow Cocone(\underline{\alpha})$. 

Using the maps $\varphi_1, \varphi_2$ with the exact triangle  $ Cocone({\alpha})\rightarrow p^*(\mathbb{T}_{{\bf X}}) \rightarrow \mathcal{O}[k]$, we obtain a triangle $ Cocone(\underline{\alpha})\rightarrow p^*(\mathbb{T}_{{\bf X}})  \rightarrow \mathcal{O}[k] \rightarrow Cocone(\underline{\alpha})[1] $ as desired.
\epf

Now, using Observation \ref{new triangle}, we get an equivalence of triangles of perfect complexes %\footnote{The first row of (\ref{equivalence of triangles of perfect complexes}) is the pullback of the triangle $\mathcal{K} \rightarrow \mathbb{T}_{{\bf X}} \rightarrow L[k]$ splitting over $\bf U.$} 
on $\spec A$
\begin{equation} \label{equivalence of triangles of perfect complexes}
	\begin{tikzpicture}
		\matrix (m) [matrix of math nodes,row sep=2.5em,column sep=1.8em,minimum width=2 em] {
		p^*(\mathcal{K})	& p^*(\mathbb{T}_{{\bf X}})  &  p^*(L)[k]  \\			
			Cocone(\underline{\alpha})  & p^*(\mathbb{T}_{{\bf X}})    & \mathcal{O}[k] ,\\
		};
		\path[-stealth]
		
		%horizontal 1st row
		(m-1-1) edge  node [right] { } (m-1-2)
		(m-1-2) edge  node [right] { } (m-1-3)
		%(m-1-4) edge  node [right] { } (m-1-5)
		%(m-1-5) edge  node [right] { } (m-1-6)
		%(m-1-6) edge  node [right] { } (m-1-7)
		%(m-1-7) edge  node [right] { } (m-1-8)
		%horizontal 2nd row
		(m-2-1) edge  node [above] {  } (m-2-2)
		(m-2-2) edge  node [above] { $ \alpha $ } (m-2-3)
		%(m-2-4) edge  node [right] { } (m-2-5)
		%(m-2-5) edge  node [right] { } (m-2-6)
		%(m-2-6) edge  node [right] { } (m-2-7)
		%	(m-2-7) edge  node [right] { } (m-2-8)
		%vertical
		(m-1-1) edge  node [right] {$ \simeq_v $ } (m-2-1)
		(m-1-2) edge  node [right] { $ id $ } (m-2-2)
		(m-1-3) edge  node [right] {$ \simeq $ } (m-2-3)
		%(m-1-6) edge  node [right] { } (m-2-6)
		%(m-1-7) edge  node [right] { } (m-2-7);
		%edge [dashed,-] (m-2-1)
		;
	\end{tikzpicture} 
\end{equation} and hence an induced isomorphism $p^*(L)\simeq\mathcal{O}.$ This map identifies  $ \mathcal{S}_{{\bf{X}}} $ with the total space of $L$.

\epf

\subsubsection{Step-2: The induced shifted symplectic structure on $ \mathcal{S}_{{\bf{X}}}$}
	Recall from  \cite{Damien2} that a morphism $\bf Y \rightarrow \widetilde{E}$ of derived stacks, with $\bf \widetilde{E}$ is the total space of $E\in QCoh({\bf X})$, consists of a morphism $f: \bf Y \rightarrow X$ together with a section $s$ of $f^*E.$ If moreover, $\bf Y$ is a derived stack equipped with a $\mathbb{G}_m$-action, then a $\mathbb{G}_m$-equivariant morphism $\bf Y \rightarrow \widetilde{E}$ is given by the pair of a $\mathbb{G}_m$-equivariant morphism $f: \bf Y \rightarrow X$ and a $\mathbb{G}_m$-equivariant section $s$ of $f^*E\{1\}.$ 
In particular, the identity map $\bf \widetilde{E} \rightarrow \widetilde{E}$ corresponds to the pair of the projection map $\pi: \bf \widetilde{E} \rightarrow X$ and a ($\mathbb{G}_m$-equivariant) section of $\pi^*E\{1\}.$ 	

Letting $\bf \widetilde{E}:= T^*[k]X$, the identity map $\bf T^*[k]X \rightarrow T^*[k]X$ is then determined by the data of the projection map $\pi_{{\bf X}}: \bf T^*[k]X \rightarrow X$ with a section of $\pi_{{\bf X}}^*\mathbb{L}_{{\bf X}}[k]$ (of  weight 1 for the fiber grading). Since we have a natural map $\pi_{{\bf X}}^*\mathbb{L}_{{\bf X}}[k] \rightarrow \mathbb{L}_{{\bf T^*[k]X}}[k]$, this section induces a $ k$-shifted 1-from $ \lambda_{\bf X} $ on $ {\bf T^*[k]X}$ called the\textit{ tautological 1-form}. Moreover, \cite{Damien2} shows that the induced closed 2-form $d_{dR}\lambda_{\bf X}$ of degree $k$ on $ {\bf T^*[k]X}$ is in fact non-degenerate, and hence gives a $k$-shifted symplectic structure.

\begin{definition}\label{defn_canonical 1-form_derived} By construction, we have the projection maps $\pi_1:  \mathcal{S}_{{\bf{X}}}\rightarrow \bf X$ and $\pi_2:  \mathcal{S}_{{\bf{X}}}\rightarrow {\bf T^*[k]X}.$ 
	We  define the \emph{canonical 1-from} $\lambda$ on $\mathcal{S}_{{\bf X}}$ to be the pullback  $ \pi_2^*\lambda_{\bf X} $ of the tautological 1-form on $ {\bf T^*[k]X}$. 
	\end{definition}

%\subsubsection*{Step-3: Shifted symplectic structure on $\mathcal{S}_{{\bf X}}$} 
Let $\mathcal{S}_{{\bf X}}, \lambda$ be as above. Then $  {\omega}:= (d_{dR}\lambda, 0, 0, \dots)$ is a $k$-shifted closed 2-form on  $\mathcal{S}_{{\bf X}}$, and hence it defines a pre-$k$-shifted symplectic structure on  $\mathcal{S}_{{\bf X}}$. Now, it remains to show that $  {\omega}$ is non-degenerate.

\begin{theorem} \label{thm_Symplectization}
	Let $\bf X$ be a (locally finitely presented) derived $\K$-scheme with a $k$-shifted contact structure $(\mathcal{K}, \kappa, L)$.
The $k$-shifted closed 2-form $ {\omega}$ described above is non-degenerate, and hence the derived stack $ \mathcal{S}_{{\bf{X}}}\rightarrow \bf X$ is $k$-shifted symplectic. %which is canonically determined by the shifted contact structure of $\bf X$ (up to quasi-isomorphism). 

	We then call the pair $(\mathcal{S}_{\bf{X}}, {\omega})$  the \emph{symplectification of $ \bf X. $}
\end{theorem}

\pf %The assertion of the theorem is local, so it is enough to prove it in a neighborhood of a point. By definition, locally on $\bf X$, where $L$ is trivialized, the perfect complex $\mathcal{K}$ in the data of the $k$-shifted contact structure on $\bf X$ can be  given as a cocone of $ \underline{\alpha} $ with $\alpha $ a locally defined $k$-shifted 1-form $\alpha: \mathbb{T}_{{\bf X}} \rightarrow \mathcal{O}_{{\bf X}}[k]$ with the property that $d_{dR}\alpha|_{\mathcal{K}}$ is non-degenerate; and thus, the triangle $\mathcal{K} \rightarrow \mathbb{T}_{{\bf X}} \rightarrow L[k]$ splits locally. Throughout the proof, we will use "refined"  neighborhoods introduced in $\S \ref{section_review of shifted symplectic strcs}$; and once the local data is specified, we will fix this $ k $-contact form.

Given $k<0$ and $p\in \bf X$, apply Theorem \ref{localmodelthm} to get a refined  neighborhood  ${\bf U}= \spec A$ of $p$ with $q \in \cspec H^0(A)$ such that $ p: \spec A \hookrightarrow \bf X $ is an open inclusion, with $q \mapsto p$, and $A$ is a minimal standard form cdga. W.l.o.g., we assume that  $L$ is trivialized on $\bf U$, and hence, over $\bf U$, the induced 1-form $ \alpha:   p^*(\mathbb{T}_{{\bf X}})\rightarrow \mathcal{O}[k]$ is such that $Cocone(\underline{\alpha}) \simeq p^*(\mathcal{K}) $ and the 2-form $d_{dR}\alpha$ is non-degenerate on $p^*\mathcal{K}$. In that case, the fiber-cofiber sequence $p^*\mathcal{K} \rightarrow p^*\mathbb{T}_{{\bf X}} \rightarrow p^*(L[k])$ splits over $\bf U.$ We fix this locally defining 1-form $\alpha$ (and $u:Cocone(\underline{\alpha})\xrightarrow{\sim}p^*\mathcal{K}$) for the rest of the proof.

Denote the pullback of $\alpha$ under the open inclusion $p$ again  by $\alpha \in p^* \mathbb{L}_{{\bf X}}[k]$. From definitions, %$Cocone(\alpha_u)\simeq \mathcal{K}$, and hence 
the triple $(p, \alpha,u)$ is an element of $\mathcal{S}_{{\bf X}}(A)$. 

Recall also that we have the distinguished triangle $p^* \mathbb{L}_{{\bf X}}\rightarrow \mathbb{L}_A \rightarrow \mathbb{L}_p$, where  $ \mathbb{L}_p $ is the relative cotangent complex, such that for refined  neighborhoods, the restriction of $\mathbb{L}_A $ to $\cspec H^0(A)$ is a free finite complex of $H^0(A)$-modules (cf. Prop. \ref{proposition_L as a complex of H^0 modules}). Moreover, we have the identification $Cocone(\underline{\alpha}) \simeq \ker \underline{\alpha}$, and  $\ker \underline{\alpha}= \ker \alpha$ in $Mod_{A}$.

From Observations \ref{observation_cocones} \& \ref{observation_equivalent contact strs},  both $\alpha$ and $f\cdot \alpha$ define equivalent contact structures, up to quasi-isomorphism, for any  $f\in H^0(A)$ (after localizing $A$ at $q$ by $f$, if necessary). 

It follows that, over ${\bf U}$, we can then identify the space $\mathcal{S}_{\bf{X}}$ locally as 
$ {\bf U} \times_{{{\bf X}}}^h \mathbb{G}_m,$\footnote{This also follows from the fact that $\mathcal{S}_{{\bf X}}$ is identified with the total space of $L$. Therefore, the corresponding homotopy fibers over $p$ are equivalent.} with natural projections, where $\mathbb{G}_m=\spec B$ is the affine group scheme, with say $B:=\spec(\K [x,x^{-1}])$. After localizing $A$ at $p$ if necessary, $\varphi \in \K [x,x^{-1}]$ acts on $f\in H^0(A)$ by $f\mapsto \varphi(f,f^{-1})$.
 
Recall  that there exists a natural equivalence  $DR(A)\otimes_{\K} DR(B) \simeq DR(A\otimes_{\K} B)$ induced by the identification
\begin{equation} \label{cotangent cmpx of tensor}
	\mathbb{L}_{A\otimes_{\K}B}\simeq(\mathbb{L}_A  \otimes_{\K} B) \oplus (A \otimes_{\K} \mathbb{L}_B). 
	\end{equation}
Notice that for $q\in {\bf U}(\K)$, with $q: \spec \K=\{\star\}\rightarrow \bf U$, and $f_q\in \mathbb{G}_m(\K)\simeq \K^{\times}$, $f_q\cdot \alpha$ is also a contact form at $q.$ Thus, on the part of the space $\mathcal{S}_{\bf{X}}$ over ${\bf U}$, we define a function $f$ with values in $\K^{\times}$. Then  the canonical 1-form in Definition \ref{defn_canonical 1-form_derived} can be locally written as%, for $f\in \mathbb{G}_m(A)$  
\begin{equation} \label{local description of can 1-form}
	\lambda=f \cdot \pi^* \alpha.
\end{equation}

\begin{remark}
This local expression (\ref{local description of can 1-form}) follows from \cite[Lemma 2.1.4]{Grataloup} that the tautological 1-form  on the shifted cotangent is universal in the sense that for any $k$-shifted 1-form $\beta$ on $\bf Y$ (viewed as $\beta: \bf Y \rightarrow T^*[k]Y $) we have $\beta^*\lambda_{{\bf Y}}=\beta$. Apply this universality for the case above ${\bf Y}:={\bf U} \times_{{{\bf X}}}^h \mathbb{G}_m$ with $\beta:=\pi^*\alpha$ and the $\mathbb{G}_m$-action to get the desired expression. We also have a non-zero factor $f$ as $\mathcal{S}_{{\bf X}}$ is identified with the total space of $L$.\footnote{That is, for (non-zero) $\pi^*\alpha, \lambda \in \mathbb{L}_{p^*\mathcal{S}_{\bf{X}}}[k]$, this identification with the total space of $L$ (i.e. the space of trivializations) implies that there is a non-zero $f\in A$ s.t. $ \lambda=f \cdot \pi^* \alpha. $}
\end{remark}

\begin{lemma} \label{proposition_symplectization}
	Locally on $\bf X$, for any locally defining 1-form $\alpha$, the $k$-shifted 2-form $\omega^0:= d_{dR} \lambda$ is  non-degenerate, and hence the sequence $\omega:=(\omega^0, 0, 0, \dots)$ defines a $k$-shifted symplectic structure on $\mathcal{S}_{\bf X}.$    
\end{lemma}
\pf[Proof of Lemma \ref{proposition_symplectization}] Note first that the non-degeneracy in sense of Definition  \ref{defn_shifted non-degeneracy} can be equivalently formulated using refined local models as follows: \begin{observation} \label{observation_nondegenarcy}
A $ k $-shifted 2-form $ \gamma  $ on $ \spec A  $  for $A$ a minimal standard form cdga is non-degenerate if and only if for any non-zero vector $v\in \mathbb{T}_A |_{q\in \cspec H^0(A)}$, there exists a non-zero vector $w\in \mathbb{T}_A |_{q\in \cspec H^0(A)}$ such that $\iota_w \iota_v \gamma \neq 0.$ For a sketch, see the note
\footnote{When restricted to $\cspec H^0(A)$, the induced morphism $ \mathbb{T}_{A} \rightarrow \Omega^1_{A}[k], \ Y \mapsto \iota_{Y} \gamma, $ in Definition  \ref{defn_shifted non-degeneracy} is just a map of finite complexes of free $ H^0(A)$-modules. And, at $q\in \cspec H^0(A)$, both $ \mathbb{T}_{A}|_q,  \Omega^1_{A}|_q$ are complexes of $\K$-vector spaces. For non-degeneracy, we require  the contraction map to be a (degree-wise) quasi-isomorphism. Recall that localizing $A$ at $ q $ if necessary, we may assume that the induces map is indeed an (degree-wise) isomorphism near $q$. Therefore, Observation \ref{observation_nondegenarcy} is just an analogous result from linear algebra.}.
\end{observation}
  Now, to prove that $\omega^0:= d_{dR} \lambda$ is  non-degenerate, we use Observation \ref{observation_nondegenarcy}. That is,
   it suffices to show that, at $q\in\bf U,$ for any non-vanishing (homogeneous) vector field $\sigma\in (\mathbb{L}_{A\otimes_{\K}B})^{\vee}$, there is  a vector field $\eta \in (\mathbb{L}_{A\otimes_{\K}B})^{\vee}$ such that $\iota_{\eta} (\iota_{\sigma} d_{dR} \lambda) \neq 0.$

For the rest of the proof,  we will also assume that the shift $k$ is odd, say $k=-2\ell -1$ for  $\ell\in \N$, to provide more explicit representations for vector fields of interest.  In fact, our constructions will be independent of the choice of the shift and the corresponding  graded variables. 

Localizing $A$ at $q$ if necessary,  choose the graded variables $x_j^{-i}, y_j^{k+i}, z^k$ on $\bf U$ as before so that $A$ is a standard form cdga over $A(0)$ freely generated by these graded variables, such that
\begin{align*}
	\ker \alpha |_{\cspec H^0(A)} &=\big \langle \partial/\partial x_1^{-i}, \dots, \partial/\partial x_{m_i}^{-i},  \partial/\partial y^{k+i}_1, \dots, \partial/\partial y^{k+i}_{m_i} : i=0,1, \dots, \ell \ \big \rangle_{H^0(A)}, \\
	Rest|_{\cspec H^0(A)} &=\big \langle \partial/\partial z^k \big \rangle_{H^0(A)} . 
\end{align*} 

%--------------to be updated--------
\begin{observation} \label{observation_summands}
	As $f|_{q\in \cspec H^0(A)} \in \K^{\times} $  and $\mathbb{L}_A=\Omega_A^1$ is a $A$-module, the first summand of Equation (\ref{cotangent cmpx of tensor}) can be equivalently written as, when restricted to ${\cspec H^0(A)}$,
	\begin{equation*}
		\mathbb{L}_A  \otimes_{\K} B|_{\cspec H^0(A)} \simeq \langle d_{dR}x_j^{-i}, d_{dR}y_j^{k+i}, d_{dR}z^k\rangle_{H^0(A)} \otimes_{\K} \K[f, f^{-1}] \simeq \langle d_{dR}x_j^{-i}, d_{dR}y_j^{k+i}, d_{dR}z^k\rangle_{H^0(A)}.
	\end{equation*} Using a cofibrant replacement of $B$ if necessary, the second summand of Equation (\ref{cotangent cmpx of tensor}) is just equivalent to $ H^0(A) \otimes_{\K} \langle d_{dR}f \rangle_{B}.  $
\end{observation}

Now, using the natural splitting \ref{splitting of complexes} and Observation \ref{observation_summands} for the complex in Equation (\ref{cotangent cmpx of tensor}), we then have, when restricted to $\cspec H^0(A)$, 
\begin{equation} \label{splitting symplectization}
	(\mathbb{L}_{A\otimes_{\K}B})^{\vee} \simeq \big(\ker \alpha  \oplus Rest \big) \oplus \big(H^0(A) \otimes_{\K}\langle \partial/\partial f \rangle_{B}\big).
\end{equation}

Using the splitting in Equation (\ref{splitting symplectization}) with  Observation \ref{observation_nondegenarcy}, we prove the statement case by case. We first note   that for any  $\eta, \sigma \in (\mathbb{L}_{A\otimes_{\K}B})^{\vee}$ at $q$, direct computations  give \begin{equation*}
\iota_{\eta} (\iota_{\sigma} d_{dR} \lambda)= \mp (d_{dR}f)(\sigma) \alpha(\pi_*\eta) \mp (d_{dR}f)(\eta) \alpha(\pi_*\sigma)  \mp f\cdot d_{dR}\alpha (\pi_* \sigma, \pi_* \eta).
\end{equation*} From Equation (\ref{splitting symplectization}), it is enough to consider the following cases:
\begin{enumerate}
	\item If $\sigma \in \ker \alpha$,  then we have %$ \pi_* \sigma \simeq \sigma $ and 
	$ \iota_{\eta} (\iota_{\sigma} d_{dR} \lambda) = \mp f\cdot d_{dR}\alpha (\sigma, \pi_* \eta)$. Since $d_{dR}\alpha |_{\ker \alpha}$ is non-degenerate by the contactness condition on $\alpha$ (and $f\in A^{\times}$), it is enough to take $\eta$ to be any non-zero vector in $\ker \alpha$.
	\item If $\sigma \in Rest,$ then we get %$ \pi_* \sigma \simeq \sigma $ and 
	$ \iota_{\eta} (\iota_{\sigma} d_{dR} \lambda) =\mp (d_{dR}f)(\eta) \alpha(\sigma). $ Here $\alpha(\sigma) \neq 0 $ since $\sigma \in Rest.$ Thus, it is enough to take $\eta$ to be any non-zero vector in $ H^0(A) \otimes_{\K}\langle \partial/\partial f \rangle_{B} $ so that $(d_{dR}f)(\eta) \neq 0.$
	\item If $\sigma \in H^0(A) \otimes_{\K}\langle \partial/\partial f \rangle_{B},$ then $ \iota_{\eta} (\iota_{\sigma} d_{dR} \lambda) = \mp (d_{dR}f)(\sigma) \alpha(\pi_*\eta).$ Observe that $ (d_{dR}f)(\sigma)\neq 0,$ so it suffices to take $\eta$ to be any non-zero vector in $ Rest $ so that $ \alpha(\pi_*\eta)\neq 0. $
\end{enumerate} This completes the proof of Lemma \ref{proposition_symplectization}, and hence that of Theorem \ref{thm_Symplectization}.

\epf

\begin{remark}
	The proofs of Lemma \ref{proposition_symplectization} and Theorem \ref{thm_Symplectization} will still be valid for the other values of $k$. %$ (a) \ k\equiv0 \mod 4, \text{ and } (b) \ k\equiv2 \mod 4  $ as well. As mentioned before, the other cases in fact involve modified versions of  the graded variables with some possible extra terms (see Equations (\ref{new local variables for k=-4l}) $ - $ (\ref{defn_phiv2})). 
	In fact, it clear to see that coordinates do not play any significant role in the proofs, rather than just providing explicit representations for the splitting. 
	
	In short, Lemma \ref{proposition_symplectization} has indeed a coordinate-free proof, and so does Theorem \ref{thm_Symplectization}. Thus, using the same terminology as before, we say that the pair $(\mathcal{S}_{\bf{X}}, \omega)$ above is the \textit{symplectification of the $k$-shifted contact derived $\K$-scheme $ \bf X $} for any $k<0.$ Note also that this result  is in fact canonical up to  quasi-isomorphism by construction.
\end{remark}

\section{Concluding remarks} \label{section_Artin stacs}

We conclude this paper with the following comments on the possible ``stacky" generalizations of the main results presented in this paper and more.% and future directions:

\begin{remark}
	It should be noted that ``stacky" generalizations of the results in \cite{Brav} are also available in the literature.  Ben-Bassat, Brav, Bussi and Joyce \cite{BenBassat etall} extend the results of \cite{Brav} from  derived schemes to the case of (locally finitely presented) derived Artin $\K$-stacks. %By a \emph{derived Artin $\K$-stack},  we essentially mean an object $\bf X$ of $ dSt_{\K}$ possessing an ``atlas" (smooth of some relative dimension) near each point of $\bf X$. That is, for each such object $\bf X$, we require the existence of a ``smooth" surjection $\varphi: U \rightarrow \bf X$ (of some relative dimension $m$), where $U$ is a derived scheme. %Here, we call such morphism \textit{an atlas.}
	In short, Ben-Bassat, Brav, Bussi and Joyce \cite{BenBassat etall} proved that derived Artin $\K$-stacks also have nice local models in some sense. Parts of the results from \cite[Theorems 2.8 \& 2.9]{BenBassat etall} in fact give the generalization of Theorem \ref{localmodelthm} to the case of derived Artin $\K$-stacks. They also proved that every shifted symplectic derived Artin $\K$-stack admits the so-called ``Darboux form atlas" \cite[Theorem 2.10]{BenBassat etall}. That is, their result extends Theorem \ref{Symplectic darboux} from  derived $\K$-schemes to  derived Artin $\K$-stacks. 
	
	Regarding the study of shifted contact structures, we note that Maglio,  Tortorella and Vitagliano \cite{Vitagliano2024} have recently introduced and studied \emph{$0$-shifted} and \emph{$+1$-shifted contact structures} on \emph{differentiable stacks}, thus providing the foundations of {shifted contact geometry} in the stacky context.%That is, near each $x\in \bf X$ one can find a ``minimal smooth atlas" $\varphi: U \rightarrow \bf X$ with $U=\spec A$ an affine derived scheme such that $(U, \varphi^*(\omega))$ is in a standard Darboux-type form.
	
	%\cite{PTVV} also defines the notions of \emph{isotropic and Lagrangian structures} on a given morphism $f: {\bf Y}\rightarrow ({\bf X}, \omega_{\bf X})$ of derived schemes with shifted symplectic target. Regarding the possible local models for these structures, Joyce and Sofronov \cite[Theorem 3.7]{JS} showed that every Lagrangian $f: {\bf Y}\rightarrow ({\bf X}, \omega_{\bf X})$ in $k$-shifted symplectic derived scheme $ ({\bf X}, \omega_{\bf X}) $, with $k<0,$ is locally modeled on explicit ``Lagrangian Darboux forms". This essentially provides a $k$-shifted version of the classical Weinstein Lagrangian neighborhood theorem.
\end{remark}	
	In the sequels,  our goals will be to extend the main results of this paper %(Theorems \ref{contact darboux} and \ref{thm_Symplectization}) 
	from derived schemes to the more general case of derived Artin stacks and to discuss more on the  theory of shifted contact derived spaces, such as \emph{Legendrians} and their local behavior, and sample constructions of contact stacks. %We begin with some basic definitions and results from \cite{BenBassat etall}. Later, we give two corollaries about the desired generalizations (cf. Corollaries \ref{cor_darboux for Artin} \& \ref{cor_symplectization for Artin}). 
	In that respect, we propose:
	\begin{claim} 
		Theorem \ref{contact darboux} and Theorem \ref{thm_Symplectization} still hold for  (locally finitely presented) $k$-shifted contact derived Artin $\K$-stacks  with $k<0$.%  with a $k$-shifted contact structure for $k<0$.  
	\end{claim}
\paragraph{The foundational trilogy.} This is the first in a series of papers \cite{kib2,kib3}. There are other ongoing projects and more to come... In brief, the second paper \cite{kib2} proves the claim above and gives several non-trivial constructions of contact derived stacks. In \cite{kib3},  we formally introduce \emph{Legendrians} in the derived context and prove a \emph{Legendrian-Darboux-type} theorem for the Legendrians in contact derived schemes.

\paragraph{Addendum.} In the latest version \bfem{(Feb 2025)} of this document, Example 3.10 and the proof of Theorem 3.13 have been revisited (along with minor edits).  The need for additional corrections and clarifications has emerged during the reviewing process of \cite{kib2}. 
Similar modifications for  \cite{kib2} will appear in  its final version. For details, see the addendum on the author's webpage. 

I wish to warmly thank the anonymous referee(s) for their valuable comments and suggestions, which helped a lot and improved the quality of the original manuscript(s). 
%\begin{conjecture} %Legendrians in derived contact stacks
%	There is a well-defined description of a \textit{Legendrian structure} on a morphism $f: \bf Y \rightarrow X$ of derived Artin stacks, with a derived contact target,  such that there also exists a Legendrian-Darboux neighborhood theorem for such structures.
%\end{conjecture}

%\section*{Acknowledgments}
%I would like to thank Ali Ulaş Özgür Kişisel and Bayram Tekin for their useful comments and suggestions, which are helpful to clarify and improve ideas presented in this paper.
%\newpage

%-----------------------------------------------------

%\bibliography{xbib}

\end{document}